\documentclass[11pt, oneside]{amsart}   	

\usepackage{geometry}
\usepackage{enumitem}
\usepackage{amsmath, amssymb}
\usepackage{mathrsfs}
\usepackage{mathabx}
\usepackage{amsthm}
\usepackage{epsfig}
\usepackage{epic,eepic} 
\usepackage[T1]{fontenc}
\usepackage[utf8]{inputenc}
\usepackage[center]{caption}
\usepackage{graphicx}				
\usepackage{fancyhdr}
\usepackage{todonotes} 	

\newtheorem{thm}{Theorem}[section]

\newtheorem{ques}[thm]{Question}

\newtheorem{conj}[thm]{Conjecture}
\newtheorem{prin}[thm]{Principle}

\newtheorem{rem}[thm]{Remark}

\theoremstyle{definition}

\newtheorem{eg}[thm]{Example}

\numberwithin{equation}{section}

\title{The Brascamp-Lieb inequality and its influence on Fourier analysis}
\author{Ruixiang Zhang}
\begin{document}

\maketitle

\begin{center}
\emph{Dedicated to Elliott Lieb on the occasion of his 90th birthday.}
\end{center}

\begin{abstract}
    Brascamp-Lieb inequalities have been important in analysis, mathematical physics and neighboring areas. Recently, these inequalities have had a deep influence on Fourier analysis and, in particular, on Fourier restriction theory. In this article we motivate and explain this connection. A lot of our examples are taken from a rapidly developing subarea called ``decoupling''. It is the author's hope that this article will be accessible to graduate students in fields broadly related to analysis.
\end{abstract}

Motivated by problems in both mathematics and physics, Elliott Lieb and his collaborators made deep contributions to a family of inequalities known as the Brascamp-Lieb inequalities today. An important special case of these inequalities is the Loomis-Whitney inequality which has a very intuitive geometric meaning. In loose terms, this inequality tells us, for example, if a set in $\Bbb{R}^3$ satisfies that all of its projections onto coordinate $x_i x_j$-planes ($1 \leq i < j\leq 3$) have area at most $100$, then the set should have volume at most $1000$. Brascamp-Lieb inequalities describe similar phenomena with a general set of projections in a general $\Bbb{R}^n$ and we have very good understanding of them now thanks to the work of Lieb \cite{lieb1990gaussian} and subsequent works.

Interestingly, variants of Brascamp-Lieb inequalities allowing ``perturbations'' attracted a lot of attention of Fourier analysts in the past 15 years. The proofs of these variants by Bennett-Carbery-Tao \cite{bennett2006multilinear} and Bennett-Bez-Flock-Lee \cite{bennett2018stability} are built both on Lieb and others' work on the ``traditional'' Brascamp-Lieb, as well as new ideas. They are considered major achievements themselves, and are also strongly motivated by problems from Fourier analysis. It turns out that these variants fundamentally revolutionized people's ways of approaching a class of problems in Fourier analysis, which in turn have deep connections to applications in number theory, combinatorics, geometric measure theory and dispersive partial differential equations. Much progress has been made due to this new tool,  known as multilinear Kakeya inequalities 
or perturbed Brascamp-Lieb inequalities in general.

In this article, we will introduce Brascamp-Lieb inequalities and their perturbed counterparts. Next we will introduce a class of problems in Fourier analysis known as ``Fourier restriction type problems'' as our main interest. Via a very intuitive geometric observation, we will then be able to see the strength of (perturbed) Brascamp-Lieb inequalities in Fourier restriction type problems, highlighting their impacts on a powerful estimate called ``decoupling''. Finally, we will explain a very recent generalization  in \cite{maldague2019regularized} of the Brascamp-Lieb inequalities and its new applications in decoupling.

As a final remark, many theorems we introduce will have their ``honest versions'' that require some quite complicated setup, and at the same time also have ``morally true'' versions that are much less technical and often serve as intuitions used by people in their research. In this introductory article, we often choose to only state the ``morally true'' version and mention references where the rigorous version can be found. We hope that graduate students and other general audience in analysis can benefit more from this choice of presentation style.

\section{From Loomis-Whitney to Brascamp-Lieb}\label{Introduction}

The Loomis-Whitney inequality  concerns the relation between the measure of a set in $\Bbb{R}^n$ and the measures of all of its projections to coordinate hyperplanes. Let $S \subset \Bbb{R}^n$ be a compact set.\footnote{We let $S$ to be compact to avoid measurability issues.} We denote the orthogonal projection to the $x_1\cdots \hat{x}_j \cdots x_n$-coordinate hyperplane by $P_j$ and use $|\cdot|$ to denote the Lebesgue measure for a set in some Euclidean space. Then we have
\begin{thm}[\cite{loomis1949inequality}]
\begin{equation}\label{LW}
    |S| \leq |P_1 (S)|^{\frac{1}{n-1}}|P_2 (S)|^{\frac{1}{n-1}}\cdots |P_n (S)|^{\frac{1}{n-1}}.
\end{equation}
\end{thm}

In two dimensions, (\ref{LW}) is immediate from the observation that $S \subset P_1(S) \times P_2(S)$ and Fubini. For higher $n$, this inequality becomes less trivial. The proof is still elementary and is a very good exercise for undergraduate analysis courses.

(\ref{LW}) has a very natural ``weighted'' generalization to (\ref{gLW}) below, often referred to as the Loomis-Whitney inequality. It can be proved in the same way as (\ref{LW}) and is often more useful. For non-negative functions $f_1, \ldots, f_n$ in $(n-1)$ variables, one has

\begin{equation}\label{gLW}
    \int_{\Bbb{R}^n} \prod_{j=1}^n f_j (x_1, \ldots, \hat{x}_j, \ldots, x_n)^{\frac{1}{n-1}} \leq \prod_{j=1}^n \left(\int_{\Bbb{R}^{n-1}} f_j\right)^{\frac{1}{n-1}}.
\end{equation}

Taking $f_j (x_1, \ldots, \hat{x}_j, \ldots, x_n)$ as $1_{P_j (S)}$, we see (\ref{LW}) follows from (\ref{gLW}).

Since (\ref{LW}) and (\ref{gLW}) seem quite interesting in their own rights, it is natural to ask whether these continue to hold for other choices of positive exponents on the right hand side or for other families of (projections, exponents) in general. Let us for example think about modifying the tuple of exponents $(\frac{1}{n-1}, \ldots, \frac{1}{n-1})$ on the right hand side of (\ref{LW}). If we play with examples a bit, we realize that if we rescale $S$ by $\lambda$ then $|S|$ is rescaled by $\lambda^n$ and each $|P_j (S)|$ is rescaled by $\lambda^{n-1}$. Hence all the exponents on the right hand side must add up to $\frac{n}{n-1}$. Moreover, one realizes that product sets $S= S_1 \times \cdots \times S_n$ saturates the bound in (\ref{LW}). Since the measures of $S_j$ can be arbitrary, we see no other choices of the exponents on the right hand side of (\ref{LW}) will make it hold, even when we replace ``$\leq$'' by ``$\lesssim$''.\footnote{This notation means that the left hand side is allowed to be no greater than a constant multiple of the right hand side.} A similar discussion applies to (\ref{gLW}).

The question about other projections is more difficult. A first observation is: If we modify $P_1, \ldots, P_n$ to be arbitrary projections to $(n-1)$ dimensional subspaces, then as long as their kernels are linearly independent, we know that (\ref{LW}) and (\ref{gLW}) continue to hold, with ``$\leq$'' replaced by ``$\lesssim$'', by a linear transform. As a more ambitious question, what can be said for an arbitrary set of projections? By considering projections of other forms (e.g. the ``Fubini'' example of $n$ projections to $n$ coordinate axes) one realizes the need of adjustments to the exponents to get a valid inequality (with $\lesssim$). And by looking at very degenerate examples (e.g. with only one nontrivial projection on the right hand side) we see that given the set of projections, it is entirely possible that no choice of exponents will yield a valid inequality (with $\lesssim$).

This motivates us to ask two questions about possible generalizations of the Loomis-Whitney inequalities. 
Assuming we are given an arbitrary set or orthogonal projections $B_j (1 \leq j \leq m)$ from $\Bbb{R}^n$ to subspaces $H_j$ and exponents $p_1, \ldots, p_m > 0$, then
\begin{ques}\label{BLques1}
Let $f_1 , \ldots, f_m$ be non-negative functions. For which $m \in \Bbb{Z}^+, p_1, \ldots, p_m >0$ and orthogonal projections $B_1, \ldots, B_m$ from $\Bbb{R}^n$ to subspaces $H_1, \ldots, H_m$ is there a finite constant $C$ such that
\begin{equation}\label{BL}
    \int_{\Bbb{R}^n} \prod_{j=1}^m f_j (B_j x)^{p_j} \leq C \prod_{j=1}^m \left(\int_{H_j} f_j\right)^{p_j}
\end{equation}
holds?
\end{ques}
\begin{ques}\label{BLques2}
If (\ref{BL}) holds, what is the sharp value of $C$ in it?
\end{ques}

Besides Loomis-Whitney, other important inequalities such as H\"{o}lder's inequality and Young's convolution inequality are also special cases of (\ref{BL}). We leave these connections to the readers as good exercises. Historically, a lot of early studies on Questions \ref{BLques1}, \ref{BLques2} and the related extremizer problem 
are motivated by generalizations of Young's inequality and problems in physics (for example see Section 5.2 of \cite{brascamp1976best}). In particular, a fair number of study was on the  ``rank one'' case, i.e. when all $\dim H_j = 1$. See e.g. Theorem 1 in \cite{brascamp1976best}.

A  fundamental contribution was made by Lieb. In \cite{lieb1990gaussian} he proved:

\begin{thm}[Lieb \cite{lieb1990gaussian}]\label{Liebthm}
To answer questions \ref{BLques1} and \ref{BLques2}, it suffices to test the inequality (\ref{BL}) with all $f_j$ being \emph{centered Gaussians} $e^{-A_j(x)}$ where $A_j$ is a positive definite quadratic form on $H_j$. In other words, the supremum of the ratio between the left hand side of (\ref{BL}) and the right hand side (allowed to be infinity) 
for all choices of $f_j$ coincides with the same supremum with each $f_j$ further restricted to be a centered Gaussian.
\end{thm}

Lieb's theorem \ref{Liebthm} is a breakthrough in the sense that it illustrates the special role of Gaussians in Questions  \ref{BLques1} and \ref{BLques2}, and reduces those questions to finite dimensional optimization problems. As a remark, even though in the previously discussed Loomis-Whitney setting it may not be so natural to expect Gaussians to play a big role, it had been known that they play central roles in the setting of Young's convolution inequality. See Babenko \cite{babenko1961inequality}, Beckner \cite{beckner1975inequalities} and Brascamp-Lieb \cite{brascamp1976best}.

Lieb actually proved several more general theorems concerning convolutions or multilinear forms with Gaussian kernels in \cite{lieb1990gaussian}. Roughly speaking, his clever proof (for Theorem \ref{Liebthm}) first establishes existence of extremizers and then notice that by taking tensor products of lower dimensional extremizers one gets a higher dimensional extremizer. The higher dimensional problem has an $O(2)$ symmetry so rotations of this extremizer are again an extremizer. But by carefully applying Minkowski's inequality, all higher dimensional extremizers must have product structure and it is impossible for rotations of products of non-Gaussian extremizers to again have product structure. In more general cases a crucial step is iterating of this ``tensor product + rotation'' construction and using an argument related to the central limit theorem.

As mentioned before, there has been a lot of study on special cases of Theorem \ref{Liebthm}, sometimes with other related purposes or different approaches. See e.g. Beckner \cite{beckner1975inequalities}, Brascamp-Lieb \cite{brascamp1976best}, Ball \cite{ball1989volumes}, Finner \cite{finner1992generalization} and Barthe \cite{barthe1998reverse}. Nevertheless, the most general Theorem \ref{Liebthm} provides a powerful unified result in all cases and is a monumental progress.

In 2005, Bennett-Carbery-Christ-Tao proved the following theorem that completely answers Question \ref{BLques1}.

\begin{thm}[Bennett-Carbery-Christ-Tao \cite{bennett2008brascamp}, see also \cite{bennett2005finite}]\label{BCCTfinite}
Given $(B_1, \ldots, B_m, p_1, \ldots, p_m)$. (\ref{BL}) holds if and only if we have both the \emph{scaling condition}:
\[n = \sum_{j=1}^m p_j \dim H_j\] and the \emph{dimension condition}: \[\dim V \leq \sum_{j=1}^m p_j \dim (B_j V), \forall \text{ subspace } V \subset \Bbb{R}^n.\]
\end{thm}

The necessity of the two conditions is easy to see. For example the scaling condition can be seen in the same way as 
before with the Loomis-Whitney case. Remarkably, Theorem \ref{BCCTfinite} states that these two conditions together are also sufficient for (\ref{BL}) to hold.

Theorem \ref{BCCTfinite} was proved by a heat flow method. To put this approach under some context, in earlier time people used rearrangement inequalities (cf. Brascamp-Lieb-Luttinger \cite{brascamp1974general}) or mass transfer inequalities (used in a related study of Barthe \cite{barthe1998reverse}) to prove special cases of  (\ref{BL}). The core idea of these methods is deforming arbitrary functions to ``good'' ones (such as Gaussians), proving some favorable monotonicity results along the deformation to reduce the problem to that of ``good'' functions, and in the endgame doing a round of computation involving good functions. A new method of this kind using heat flows was first used by Carlen-Lieb-Loss \cite{carlen2004sharp} to obtain a sharp Young's inequality on the sphere. Bennett-Carbery-Christ-Tao also used heat flow methods to prove their Theorem \ref{BCCTfinite}, which (though rediscovered independently) can be viewed as a great extension of the arguments in \cite{carlen2004sharp} to the full general higher-rank case.

We would like to mention a third important work stating that inequality (\ref{BL}) is stable under perturbation. The main theorem there is informally stated below.

\begin{thm}[Bennett-Bez-Flock-Lee \cite{bennett2018stability}]\label{BLstab}
If (\ref{BL}) holds for some $(B_1, \ldots, B_m, p_1, \ldots, p_m)$, then it also holds with a uniformly bounded implied constant for $(\tilde{B}_1, \ldots, \tilde{B}_m, p_1, \ldots, p_m)$, where each $\tilde{B}_j$ is a sufficiently small perturbation of $B_j$, respectively. Here the smallness may depend on $(B_1, \ldots, B_m, p_1, \ldots, p_m)$.
\end{thm}

The proof of Theorem \ref{BLstab} is a further development based on Lieb and Bennett-Carbery-Christ-Tao's theorems.

Following the literature, we call (\ref{BL}) a \emph{Brascamp-Lieb inequality} when it holds. For general $(B_1, \ldots, B_m, p_1, \ldots, p_m)$, we call the sharp constant $C$ (which may or may not be infinity) such that (\ref{BL})
holds to be the \emph{Brascamp-Lieb constant} (usually denoted as $BL(\mathbf{B}, \mathbf{p})$) associated to the \emph{Brascamp-Lieb datum} $(\mathbf{B}, \mathbf{p})$ with $\mathbf{B} = (B_1, \ldots, B_m)$ and $\mathbf{p} = (p_1, \ldots, p_m)$.

\subsection{Additional historical remarks}
Only a restricted form of the Brascamp-Lieb inequality appeared in the original paper \cite{brascamp1976best} by Brascamp-Lieb and to the best of the author's knowledge the general inequality first appeared in Lieb's work \cite{lieb1990gaussian}. There are other important inequalities proved in \cite{brascamp1976best} that we will not be able to touch in this article. For example they proved a converse inequality that implies the celebrated Pr\'{e}kopa-Leindler theorem.

The Brascamp-Lieb inequality and its geometric form by Ball \cite{ball1989volumes, ball1991volume} and Barthe \cite{barthe1998reverse} are important, e.g. to prove smoothing estimates for heat kernels of (possibly magnetic) Schr\"{o}dinger operators (\cite{bourga1991besicovitch, lieb1990gaussian, loss1997optimal} to only mention three), in kinetic theory (e.g. \cite{bonetto2018entropy}) and other fields.

\section{Fourier restriction type problems}\label{FRTprob}

We now turn to another main topic of this article: Fourier restriction type problems. These problems form an important subfield of Fourier analysis. Not only are they important in their own rights but they also have a lot of connections to other mathematical subjects including analytic number theory, dispersive partial differential equations, combinatorics, spectral theory and geometric measure theory. Next we mention some examples of these problems and connections, and explain a discovery people made in the last 15 years on how Brascamp-Lieb inequalities (or their variants) can be crucial tools in partial or full solutions to Fourier restriction type problems.

Good references in Fourier restriction include for example Wolff \cite{wolff2003lectures}, Tao \cite{tao2001rotating, tao2004recent}, Stovall \cite{stovall2019waves}, Mattila \cite{mattila2015fourier} and most recently Demeter \cite{demeter2020fourier}. 

\subsection{The problems}

Let us start by recalling some basic facts about the Fourier transform. For a Schwartz function $\phi \in \mathscr{S} (\Bbb{R}^n)$, its Fourier transform $\hat{\phi}$ is defined as \[\hat{\phi} (\xi) = \int_{\Bbb{R}^n} e^{-2\pi \mathrm{i}x \cdot \xi} \phi (x) \mathrm{d}x\] and also belongs to $\mathscr{S} (\Bbb{R}^n)$. One can use the above definition to extend the Fourier transform to the space $\mathscr{S}' (\Bbb{R}^n)$ of tempered distributions. We have the Plancherel theorem: \[\|\hat{\phi}\|_{L^2 (\Bbb{R}^n)} = \|\phi\|_{L^2 (\Bbb{R}^n)}\] and the Hausdorff-Young inequality\footnote{Here $p'$ denotes the dual exponent of $p$: $\frac{1}{p}+\frac{1}{p'} = 1$.}: \[\|\hat{\phi}\|_{L^{p'} (\Bbb{R}^n)} \leq \|\phi\|_{L^p (\Bbb{R}^n)}, 1 \leq p \leq 2.\]

We can now introduce the subject of \emph{Fourier restriction}. Assume we have a function $g \in L^{\infty} (\Bbb{R}^n)$ such that $\text{supp }\hat{g} \subset S$ where $S$ is a compact submanifold\footnote{Often with boundary.}. A fundamental observation by E. M. Stein (cf. \cite{fefferman1970inequalities, fefferman1995selected}) is that when $S$ is ``curved'', very often one can hope for certain nontrivial (and often difficult) estimates on $\|g\|_{L^p}$ for some $p \geq 1$. Here the $L^p$ norm is taken with respect to either the Lebesgue measure, or a certain weight, such as a Frostman measure on a fractal set. A \emph{Fourier restriction} type problem requires one to obtain this kind of $L^p$ estimates.

To be more precise, let us look at a few important examples and some of their surprising connections to other fields. The single most important problem of this kind concerns Stein's observation of the ``Fourier restriction'' phenomenon. Take $S$ to be the truncated unit paraboloid \[P^{n-1} = \{\xi_{n} = \sum_{j=1}^{n-1} \xi_j^2, |\xi_j| \leq 1, \forall 1 \leq j \leq n-1\}.\]

From the definition, we know that if $\phi \in L^1 (\Bbb{R}^n)$, then $\hat{\phi}$ is bounded, continuous and defined everywhere on $P^{n-1}$. On the other hand, since Plancherel is an $L^2$ isometry, generally one can not hope to make sense of $\hat{\phi}$ as an integrable function on $P^{n-1}$ (that has zero measure) when we merely have $\phi \in L^2 (\Bbb{R}^n)$. Stein's remarkable discovery is that for certain $q \in (1, 2)$, whenever $\phi \in L^{q} (\Bbb{R}^n)$, one should expect $\hat{\phi}$ to makes sense as an integrable function on $P^{n-1}$ (with respect to the natural surface measure $\mathrm{d}\sigma$). This leads to  Stein's Fourier restriction conjecture \cite{stein1979some}

\begin{conj}[Fourier restriction]\label{restrconj}
For $1 \leq q < \frac{2n}{n+1}$ and every $\phi \in \mathscr{S} (\Bbb{R}^n)$, \begin{equation}\label{restrconjineq}
    \|\hat{\phi}\|_{L^1 (P^{n-1}, \mathrm{d}\sigma)} \lesssim \|\phi\|_{L^q (\Bbb{R}^n)}.
\end{equation}
\end{conj}

Conjecture \ref{restrconj}, if true, has the immediate corollary that for $1 \leq q < \frac{2n}{n+1}$ and every $\phi \in L^q (\Bbb{R}^n)$, $\hat{\phi}$ can be meaningfully restricted to $P^{n-1}$ as an integrable function with respect to $\mathrm{d}\sigma$. This explains its name ``Fourier restriction''.

By duality, Conjecture \ref{restrconj} is equivalent to  the following ``Fourier extension'' conjecture that is often friendlier to work with\footnote{As a basic fact, the Fourier transform is invertible on $\mathscr{S}'$ and we use $\check{\cdot}$ to denote its inverse.}:

\begin{conj}[Fourier extension]\label{extconj}
For $p > \frac{2n}{n-1}$ and every $f \in L^{\infty} (P^{n-1}, \mathrm{d}\sigma)$, \begin{equation}\label{Fexteq}
    \|\widecheck{f \mathrm{d}\sigma}\|_{L^p (\Bbb{R}^n)} \lesssim \|f\|_{L^{\infty}}.
\end{equation}
\end{conj}

The exponent $\frac{2n}{n-1}$ in Conjecture \ref{extconj} (as well as the $\frac{2n}{n+1}$ in Conjecture \ref{restrconj}) is optimal, as  can be seen through a simple example. If one takes $f$ to be a bump function in $(\xi_1, \ldots, \xi_{n-1})$, then by the method of stationary phase, $|\widecheck{f \mathrm{d}\sigma}|$ decays like $|x|^{-\frac{n-1}{2}}$ and this decay rate is on average the best, leading to the sharpness of $\frac{2n}{n-1}$.

Conjecture \ref{extconj} is expected to hold, but as of today it is only fully known in dimension $2$ (see \cite{fefferman1970inequalities} and the endpoint version in \cite{zygmund1974fourier}) and remains widely open in higher dimensions. It is closely related to other long-standing major open problems in Fourier analysis such as the Bochner-Riesz conjecture \cite{bochner1935summation}, and is a central and fundamental unsolved problem in modern Fourier analysis itself.

Note that $g = \widecheck{f \mathrm{d}\sigma}$ has the Fourier support in $P^{n-1}$ and thus Conjecture \ref{extconj} is a typical Fourier restriction type problem that we introduced a while ago.\footnote{In fact, this is where the name ``Fourier restriction'' comes from.} Another remark is that one easily sees no interesting estimates can exist if we replace $P^{n-1}$ by open subsets of linear subspaces. Thus one sees that the presence of ``curvature'' is a key reason why estimates like (\ref{Fexteq}) may be expected. For the ubiquitous roles played by curvatures in harmonic analysis, see \cite{littman1963fourier, greenleaf1981principal} and the references therein.

We now introduce a few other Fourier restriction type problems that play important roles in neighboring subjects.

\begin{eg}
Consider the solution $u(x, t)$ to the free wave equation
\begin{equation}\label{waveeq}
u_{tt} = \Delta_x u, x \in \Bbb{R}^n, t \in \Bbb{R}
\end{equation} with initial position $u(x, 0) = f(x)$ and initial velocity $u_t(x, 0) = 0$. We remark that one can also consider a more general setting $u (x, 0) = f_1 (x), u_t (x, 0) = f_2 (x)$, where there is a story entirely parallel to the discussion below that we will not go into.

For $p>2$, upper bounds of $L^p$ norms of the solution $u$ are of much interest since it tells us how much the solution can focus. To see this, note that we have the fixed time $L^2$ norm bounded by energy conservation. For $p>2$, the  $L^p$ bound is getting big if the solution ``focuses'' as a function.

Three decades ago, Sogge \cite{sogge1991propagation} noticed that a local averaging in time (say averaging over $t \in [1, 2]$) is expected to have a ``smoothing'' effect on solutions, leading one to much better $L^p$ estimates than the best possible fixed time estimate observed by Miyachi \cite{miyachi1980some} and Peral \cite{peral1980lp}. When $p$ is only slightly larger than $2$, it is expected that the averaged time $L^p$ estimate only loses an arbitrary small amount of derivatives. This is the following \emph{local smoothing} conjecture:

\begin{conj}\label{localsmoothingconj}[local smoothing for free waves \cite{sogge1991propagation}]
If $u$ solves (\ref{waveeq}), then for $2 \leq p \leq \frac{2n}{n-1}$, \[\|u\|_{L^p (\Bbb{R}^n \times [1, 2])} \leq C_{\varepsilon} \|g\|_{W^{\varepsilon, p} (\Bbb{R}^n)}, \forall \varepsilon > 0.\] Here $W^{\varepsilon, p}$ is a  Sobolev space. Intuitively it consists of functions in $L^p$ with ``$\varepsilon$-many $L^p$ derivatives''. For a rigorous definition of the $W^{\varepsilon, p}$ norm, we can use $\|(1+\sqrt{-\Delta})^{\varepsilon} f\|_p$.
\end{conj}

If we do the space-time Fourier transform to  (\ref{waveeq}), we obtain a nice support condition that $\text{supp }\hat{u} \subset \tilde{\Gamma}^{n}$ with $\tilde{\Gamma}^{n}$ being the unit cone \[\tilde{\Gamma}^{n} = \{(x, t): |x| = |t|\}.\]

One can perform a Littlewood-Paley decomposition to $f$, i.e. a decomposition of $f$ into pieces such that the frequency support of each piece is in a unique dyadic annulus. By doing this, $u$ will be decomposed into pieces whose Fourier supports are on the intersection of $\tilde{\Gamma}^{n}$ and various dyadic annuli. This decomposition is rather ``mild'' and the original conjecture is equivalent to the problem concerning one individual piece at a time. Doing an appropriate rescaling to that piece, we arrive at a function on $\Bbb{R}^{n+1}$ whose Fourier support is in the \emph{truncated unit cone} \[\Gamma^{n} = \{(x, t): |x| = |t|, 1 \leq |x| \leq 2\}.\] We would like to know the $L^p$ norm of this function and hence we are interested in a Fourier restriction type problem associated to the truncated unit cone in the frequency space\footnote{Meaning the space where the Fourier transform is defined on.}.
\end{eg}

\begin{eg}
Throughout this article, a natural number means zero or a positive integer. More than one and a half centuries ago, Waring's problem was raised that asks: for each positive integer $k$, what is the minimal $s$ such that every natural number can be written as a sum of $s$ $k$th-powers of naturals?  It is also of interest to ask: What is the minimal $s$ such that the above holds for every sufficiently large natural number? For example, Lagrange proved in 1770 that every natural number is a sum of four squares. On the other hand, there exist arbitrary large natural numbers that cannot be written as a sum of three squares. Indeed everything that is $\equiv 7 (\mod 8)$ is not a sum of three squares.

Hilbert \cite{hilbert1932beweis} proved that $s$ can be finite for every $k$. But it is very difficult to determine the exact $s$ and number theorists are interested in obtaining good upper bounds for $s$. The celebrated circle method developed by Hardy and Littlewood (that originated in the work of Hardy and Ramanujan) is very fruitful in estimating $s$. Via this method, the counting of solutions are connected to estimates of both pointwise values and $L^p$ norms of certain periodic functions. In applications of this method it is important to understand the asymptotic count of solutions to equations like \begin{equation}\label{Circleeqndeg4}
    x_1^4+ \cdots + x_{10}^4 = y_1^4+ \cdots + y_{10}^4, x_j, y_j \in \Bbb{Z}, |x_j|, |y_j| \leq N.
\end{equation}
In the case of (\ref{Circleeqndeg4}) we expect the number of solutions to have an upper bound $O(N^{16})$ or $O_{\varepsilon}(N^{16+\varepsilon})$ predicted by a probabilistic consideration. Vinogradov \cite{vinogradov1935new, vinogradov1947method} made significant contributions to the subject. He noticed that the above expected sharp upper bounds can be deduced from sharp upper bounds of the solutions to the system \begin{equation}\label{Vleeqndeg4}
    x_1^i+ \cdots + x_{10}^i = y_1^i+ \cdots + y_{10}^i, 1\leq i\leq 4,  x_j, y_j \in \Bbb{Z}, |x_j|, |y_j| \leq N.
\end{equation}
The expected solution count to the system (\ref{Vleeqndeg4}) is $O(N^{10})$ or $O_{\varepsilon}(N^{10+\varepsilon})$, which implies the sharp upper bound for solution count to equation (\ref{Circleeqndeg4}). But (\ref{Vleeqndeg4}) turns out to be easier to study, as we have better ``curvature'' conditions to use in a corresponding Fourier restriction type problem. To see this, note that by Plancherel, the number of solutions to (\ref{Vleeqndeg4}) is simply $\|g_{N}(\mathbf{\alpha})\|_{L^{20} (\Bbb{T}^4)}^{20}$, where $\Bbb{T}^4 = (\Bbb{R}/ \Bbb{Z})^4$ is the four dimensional torus and  \[g_{N}(\mathbf{\alpha}) = \sum_{|n| \leq N} e^{2\pi\mathrm{i} (\alpha_1 n +\alpha_2 n^2 + \alpha_3 n^3 +\alpha_4 n^4)}.\]
Now by rescaling, $g_{N} (\frac{\cdot}{N}, \frac{\cdot}{N^2}, \frac{\cdot}{N^3}, \frac{\cdot}{N^4})$ can also be viewed as a periodic function on $\Bbb{R}^n$ (here $n=4$) whose Fourier support is on the \emph{truncated moment curve} \begin{equation}
    \mathcal{M}_n = \{(x, x^2, \ldots, x^n): |x| \leq 1\}.
\end{equation}
Since we care about the $L^{20}$ norm of $g_N$ on a box (when we view $g_N$ as a periodic function), we are looking at a typical Fourier restriction type problem. Since the moment curve is curved in a very non-degenerate way, one may imagine that this helps us to get good estimates for $g_N$.
\end{eg}

\begin{eg}
A very influential problem in geometric measure theory is Falconer's distance conjecture:
\begin{conj}[Falconer's distance conjecture \cite{falconer1985hausdorff}]\label{Falconerconj}
For a set $E$ in $\Bbb{R}^n$, define \[\Delta (E) = \{|x-x'|: x, x' \in E\}\subset \Bbb{R}\] to be the set of distances it generates. Then when $n>2$ and $E$ has Hausdorff dimension $>\frac{n}{2}$, we must have $|\Delta (E)| > 0$.
\end{conj}

This conjecture is widely open even in two dimensions. Because of the spherical symmetry of the distance function, it is known (Mattila \cite{mattila1985hausdorff, mattila1987spherical}, see also Liu \cite{liu20192}) that we can make progress on this problem\footnote{Which means lowering the dimension needed to make $|\Delta (E)| > 0$.} if we know good $L^2$ estimates of functions whose Fourier supports are on the unit sphere $S^{n-1}$, with respect to rescaled Frostman measure on the set. This is a very nontrivial Fourier restriction type problem.

The best known result on this conjecture is obtained by Guth-Iosevich-Ou-Wang \cite{guth2020falconer} in two dimensions, and by the author and collaborators (Du, Guth, Iosevich, Ou, Wang, Wilson) \cite{du2021weighted, du2019sharp, du2021improved} in higher dimensions. All these results as well as earlier results by Wolff \cite{wolff1999decay} and Erdo\u{g}an \cite{erdogan2005bilinear} make use of the above framework and rely on new Fourier restriction type estimates.
\end{eg}

\subsection{Basic tools: The uncertainty principle and the wave packet decomposition}

Having seen the importance of Fourier restriction problems, we now come to methods to deal with them. We will introduce a few consequences of the \emph{uncertainty principle}, which are now very standard tools. It will then become evident that Brascamp-Lieb type inequalities can provide essential assistance in the picture.

In a Fourier restriction type problem, we have a function $g$ whose Fourier support is on a submanifold $S$. As analysts, we would like to decompose $g$ into pieces. We would like to better understand each piece and understand how they add up together. A very useful tool in both the heuristic and the rigorous sides is the \emph{uncertainty principle}. This principle may have different meaning depending on contexts and we will loosely state the most relevant version to us:

\begin{thm}[Uncertainty principle, loosely stated]\label{UCTP}
If $\text{supp } \hat{g}$ is in a rectangular box $B$, then $|g|$ is essentially a constant on every translation of the dual box $B^*$ to $B$.
\end{thm}

By definition, the dual box $B^*$ is centered at $0$ and has edges parallel to corresponding edges of $B$ but with lengths equal to the inverse of the counterparts of $B$.

A rigorous versions of Theorem \ref{UCTP} can be found as Proposition 5.5 in \cite{wolff2003lectures}. The readers will see that it involves quickly decaying weights for technical reason. Its proof is not difficult: just notice $g = g*\psi_{B^*}$ where $\psi_{B^*}$ is a mollifier adapted to $B^*$.

As an example, by this principle we see that when $S$ is compact, we can pretend that $|g|$ behaves like a constant function in every unit cube.

To understand $\|g\|_{L^p}$ on the whole space or with respect to some weight in our Fourier restriction problem, it is often enough to understand pointwise or weighted bounds of this function inside a large ball of radius $R>1$ sufficiently well, and then let $R \to \infty$. For example in the context of the Fourier restriction conjecture, Tao \cite{tao1999bochner} proved a rigorous version of this reduction (see also Section 2.1 of \cite{tao2004recent} or Section 1.6 of \cite{demeter2020fourier}).

When we restrict our attention inside a ball $B_R$ of radius $R$, there is a \emph{wave packet decomposition} to decompose $g$ into simpler pieces (wave packets). Each wave packet is morally localized both in the physical space and in the frequency space at the scale that is allowed by the uncertainty principle. The idea of localizing the function in both the physical and frequency spaces goes back at least to Fefferman \cite{fefferman1973note} and C\'ordoba \cite{cordoba1977kakeya, cordoba1979maximal}. Bourgain \cite{bourga1991besicovitch} and subsequent authors systematically use this method to work on Fourier restriction type problems.

We loosely state this decomposition in the next theorem in the important case where $S = P^{n-1}$ (The case $S = S^{n-1}$ is entirely similar, which is useful in the study of e.g. Falconer's distance Conjecture \ref{Falconerconj}). As with the uncertainty principle, rigorous statements of this decomposition often involve rapidly decaying tails and will not be stated here. See e.g. Chapter 2 of \cite{demeter2020fourier} or Section 3 of \cite{guth2018restriction} for rigorous statements.

\begin{thm}[Wave packet decomposition, loosely stated]\label{WPD}
Decompose $P^{n-1}$ into $\sim R^{\frac{n-1}{2}}$ disjoint pieces $\theta$ with diameter $R^{-\frac{1}{2}}$ (called \emph{$R^{-\frac{1}{2}}$-caps}). Simple differential geometry shows each piece is contained in a rectangular box $\overline{\theta}$ of dimensions roughly $R^{-\frac{1}{2}}\times \cdots \times  R^{-\frac{1}{2}}\times R^{-1}$. For $\text{supp } \hat{g} \subset P^{n-1}$, decompose \[g = \sum_{\theta} g_{\theta}\] where each $\hat{g}_{\theta}$ is defined to be $\hat{g}$ restricted to $\theta$. For each $\theta$, tile $B_R$ by translations $T$ of $\overline{\theta}^*$ (and write $T \sim \theta$). These are rectangular boxes of dimensions $R^{\frac{1}{2}}\times \cdots \times  R^{\frac{1}{2}}\times R$ with edges parallel to the corresponding edges of $\overline{\theta}$. Then we can further decompose each \[g_{\theta} = \sum_{T \sim \theta} g_{\theta, T}\] such that:

(i) Each $|g_{\theta, T}|$ is essentially supported on $T$ and is like a constant on $T$.

(ii) Each $\hat{g}_{\theta, T}$ is roughly supported on $\theta$.

(iii) Different $g_{\theta, T}$ are morally orthogonal in the $L^2$ sense on  every $R^{\frac{1}{2}}$-ball. In other words, for a collection $I$ of $(\theta, T)$ and a ball $B_{R^{\frac{1}{2}}}$, it is approximately true that \[``\|\sum_{(\theta, T) \in I} g_{\theta, T} \|_{L^2 (B_{R^{\frac{1}{2}}})}^2 = \sum_{(\theta, T) \in I} \| g_{\theta, T} \|_{L^2 (B_{R^{\frac{1}{2}}})}^2 ."\]
\end{thm}

In the above decomposition, each $g_{\theta, T}$ is called a \emph{wave packet}. They are simple functions to understand: They have morally constant magnitude on $T$ and decay rapidly out of it. Still, we need to put them together to recover the original $g$. When wave packets overlap, we have local $L^2$ orthogonality between them but this is usually not strong enough in practice to get $L^p$ estimates for $p > 2$.  For this reason, it is in general more difficult to obtain good or sharp $L^p$ estimates for $g$ when the wave packets overlap a lot.

Fortunately, there are good reasons to believe that wave packets intersect in a limited way (known as Kakeya type phenomenon), which is purely a combinatorial statement about how tubes can overlap. If one can prove quantitative versions of this ``limited intersection'' statement, they can possibly make significant progress on Fourier restriction type estimates. Pioneering works along these lines in the setting of the Fourier restriction conjecture was done by Bourgain \cite{bourga1991besicovitch} and subsequent authors (see e.g. \cite{wolff2001sharp, tao2003sharp}). Kakeya type problems and their connections to Fourier restriction will be briefly introduced in \S \ref{reallysharpsec}, and can be found in many survey articles such as \cite{tao2001rotating}.

About 15 years ago, people started to notice that generalizations of the Loomis-Whitney inequality or the more general Brascamp-Lieb inequality are plausible, and can lead to new ``limited intersection'' statements for tubes that help in Fourier restriction theorems. This will be the main topic of the next section.

\subsection{Additional historical remarks} Conjectures \ref{restrconj} and \ref{extconj} have other equivalent forms, see e.g. Conjectures 19.5 and 19.7 in Mattila's book \cite{mattila2015fourier}.

Like the Fourier restriction conjecture, the Bochner-Riesz conjeture is also solved in dimension two (cf. \cite{carleson1972oscillatory}) but widely open in higher dimensions. For a good introduction to the Bochner-Riesz conjeture and its connection to Conjectures \ref{restrconj} and \ref{extconj}, see e.g. Lecture 5 in \cite{tao2004recent} or table 2 in \cite{tao1999bochner}.

The genesis of the circle method is in the a paper by Hardy and Ramanujan in 1917 \cite{ramanujan2015collected} on the asymptotic number of the total number of partitions of $n$. The method was then developed in great details in a series of papers by Hardy and Littlewood on ``Partitio Numerorum''. Among the series, \cite{hardy1920some, hardy1921some, hardy1922some, hardy1925some, hardy1928some} are mainly on Waring's problem. Davenport's book \cite{davenport2005analytic} has a good introduction to the circle method. 

\section{Perturbed Brascamp-Lieb inequalities}

People have found it useful to develop perturbed versions of Loomis-Whitney inequalities and more general Brascamp-Lieb inequalities in the study of Fourier restriction problems introduced in the last section. In this section we explain this connection. 

\subsection{Bilinearizing and multilinearizing the Fourier extension problem}\label{MRsubsec}
There are variants of the Fourier extension problem (Conjecture \ref{extconj}). One can ask the same question with other Lebesgue spaces on both sides of (\ref{Fexteq}). For an integrable function $f$ on $P^{n-1}$, for simplicity we define \[Ef = \widecheck{f\mathrm{d}\sigma}.\] Note that by the uncertainty principle $|Ef|$ is like a constant on every unit ball, and thus (\ref{Fexteq}) becomes stronger if one decreases the Lebesgue exponent on the left hand side.

An important variant is when one replaces $L^{\infty}$ by $L^2$ on the right hand side. It turns out that the sharp Lebesgue exponent on the left hand side has to be increased to $\frac{2(n+1)}{n-1}$. This is then the celebrated Stein-Tomas theorem which is known to be true even at the endpoint $p= \frac{2(n+1)}{n-1}$:
\begin{thm}[Stein-Tomas, \cite{tomas1975restriction, stein1986oscillatory}]\label{STthm}
\begin{equation}\label{FextL2eq}
    \|Ef\|_{L^{\frac{2(n+1)}{n-1}} (\Bbb{R}^n)} \lesssim \|f\|_{L^{2} (\mathrm{d}\sigma)}.
\end{equation}
\end{thm}

To see the sharpness of $\frac{2(n+1)}{n-1}$ in (\ref{FextL2eq}), we consider an important example known as the \emph{Knapp example}. For a large $R>1$, take $f$ to be $1$ on one $R^{-\frac{1}{2}}$-cap and zero elsewhere. Then since a cap is contained in a rectangular box of size $\sim R^{-\frac{1}{2}} \times \cdots \times R^{-\frac{1}{2}} \times R^{-1}$, an elementary computation (that agrees with the uncertainty principle heuristic) shows that $|Ef| \sim R^{-\frac{n-1}{2}}$ on a tube of length $\sim R$ and radius $\sim R^{\frac{1}{2}}$.\footnote{We sometimes work with tubes rather than comparable rectangular boxes for notational convenience.} We then see that for (\ref{FextL2eq}) to hold, the $\frac{2(n+1)}{n-1}$ cannot be dropped. In the context of the wave packet decomposition, the Knapp example 
can morally be viewed as there being only one significant wave packet while every other wave packet is negligible.

With the wave packet decomposition in mind, it is conceivable that Fourier restriction estimates with $L^2$ on the right hand side are easier to understand than those with other $L^p (p>2)$ on the right hand side. This is because in the $L^2$ world we have Plancherel and local $L^2$ orthogonality of wave packets ((iii) in Theorem \ref{WPD}). Indeed, we now have the very sharp Theorem \ref{STthm} while Conjecture \ref{extconj} has resisted many attacks and remains open today for $n \geq 3$.

Do we learn anything about Conjecture \ref{extconj} from the $L^2$-based Fourier restriction estimate (\ref{FextL2eq})? Definitely: By H\"{o}lder we know Conjecture \ref{extconj} is true with the $\frac{2n}{n-1}$ replaced by $\frac{2(n+1)}{n-1}$. Unfortunately this is the best we can get out of (\ref{FextL2eq}), since the Knapp example prevents us from further dropping the Lebesgue exponent on the left hand side.

Is there a different estimate with $L^2$ on the right hand side and some Lebesgue exponent below $\frac{2(n+1)}{n-1}$ on the left hand side that can be useful in Fourier restriction? Surprisingly, the answer is yes and the key is to ``bilinearize'' the estimate. Here we state such a bilinear restriction estimate.

\begin{thm}[Bilinear restriction \cite{tao2003sharp}]\label{bilinearrestrthm}
Fix open sets $S_1, S_2 \subset P^{n-1}$ with $\text{dist } (S_1, S_2) \gtrsim 1$. Then for $\text{supp }f_1 \subset S_1$ and $\text{supp }f_2 \subset S_2$,
\begin{equation}\label{bilinearineq}
    \||Ef_1\cdot Ef_2|^{\frac{1}{2}}\|_{L^{\frac{2(n+2)}{n}+\varepsilon} (\Bbb{R}^n)} \lesssim_{\varepsilon} \|f_1\|_{L^{2} (\mathrm{d}\sigma)}^{\frac{1}{2}}\|f_2\|_{L^{2} (\mathrm{d}\sigma)}^{\frac{1}{2}}.
\end{equation}
\end{thm}

The formulation (\ref{bilinearineq}) successfully dropped the Lebesgue exponent on the left hand side, which is expected from examples such as the following: If we take one wave packet of $f_1$ and one of $f_2$ then the intersection of their supports will essentially be a $R^{\frac{1}{2}}$-ball, much smaller than each individual support of the wave packet. This suggests that it is possible for the bilinear estimate (\ref{bilinearineq}) to have a lower exponent on the left hand side than Stein-Tomas. Indeed, Tao's proof makes essential use of the geometric observation that any such pair of wave packets can only have supports intersect roughly at an $R^{\frac{1}{2}}$-ball.

It is known that estimates like (\ref{bilinearineq}) directly imply progress on the Fourier extension Conjecture \ref{extconj} via a Whitney decomposition argument (see e.g. Chapters 3-4 of \cite{demeter2020fourier}). In particular, Theorem \ref{bilinearrestrthm} implies (\ref{Fexteq}) with the exponent on the left hand side replaced by $\frac{2(n+2)}{n}+\varepsilon$, and hence for example implies Conjectures \ref{restrconj} and \ref{extconj} for $n=2$.

Being a powerful method itself, the bilinear approach also has limits. It is a very good exercise that the exponent ${\frac{2(n+2)}{n}}$ in Theorem \ref{bilinearrestrthm} is sharp, seen by constructing an example of $f_1$ and $f_2$ using parallel wave packets tiling an $R^{\frac{1}{2}} \times \cdots R^{\frac{1}{2}} \times R \times R$-box. If one wants to continue lowering the Lebesgue exponent on the left hand side of (\ref{Fexteq}), or to prove the whole  Conjecture \ref{extconj}, new insights are needed.

What if we further multilinearize the left hand side? By doing this, we can potentially have more transversality between wave packets. With a careful setup, the  example in the last paragraph will be out of the picture and indeed the exponent on the left hand side may be further dropped! In the $n$-linear setting, the exponent can be lowered  all the way down to $\frac{2n}{n-1}$, matching the left hand side of (\ref{Fexteq}). This is the following celebrated multilinear restriction theorem by Bennett-Carbery-Tao:

\begin{thm}[Multilinear restriction \cite{bennett2006multilinear}]\label{MRthm}
Fix open sets $S_1, S_2, \ldots, S_n \subset P^{n-1}$ satisfying the \emph{transverse condition}: For any hyperplane $b$, there is $1 \leq j \leq n$ such that the normal vector of every point in $S_j$ has an angle $\gtrsim 1$ against $b$. Then for $\text{supp }f_j \subset S_j, 1 \leq j \leq n$ and every ball $B_R$ of radius $R>1$,
\begin{equation}\label{multilinearineq}
    \||\prod_{j=1}^n Ef_j|^{\frac{1}{n}}\|_{L^{\frac{2n}{n-1}} (B_R)} \lesssim_{\varepsilon} R^{\varepsilon} \prod_{j=1}^n\|f_j\|_{L^{2} (\mathrm{d}\sigma)}^{\frac{1}{n}}, \forall \varepsilon>0.
\end{equation}
\end{thm}

The multilinear restriction theorem actually holds in much more general setting: One can start with arbitrary smooth hypersurfaces $S_1, \ldots, S_n$, not necessarily on $P^{n-1}$ nor need to have any curvature individually. As long as they satisfy the transverse condition, there is an entirely analogous inequality to (\ref{multilinearineq}). A technical remark is that here the restriction to $B_R$ and the loss of an $R^{\varepsilon}$ factor are no big losses, thanks to epsilon-removal lemmas (more on them in \S \ref{AHR3}).

It is exciting that we are allowed to have the magical exponent $\frac{2n}{n-1}$ on the left hand side of (\ref{multilinearineq}) to match the sharp exponent in Conjecture \ref{extconj}. However unlike the bilinear case, Theorem \ref{MRthm} does not immediately imply Conjecture \ref{extconj} or any partial progress on it. It is generally believed that Theorem \ref{MRthm} is much weaker than Conjecture \ref{extconj}.

Nevertheless, Theorem \ref{MRthm} is surely interesting in its own right and it felt to people it should result in new  progress on Conjecture \ref{extconj}. Indeed, a few years later, improvements on Conjecture \ref{extconj} from Theorem \ref{MRthm} were finally made possible \cite{bourgain2011bounds, guth2018restriction}. As of today, Theorem \ref{MRthm} has helped people make tremendous progress on, or fully resolve numerous other Fourier restriction type problems, and is one of the most influential tools in Fourier restriction. This will be the topic of the next few sections.

The single most important tool in the proof of Theorem \ref{MRthm} is a perturbed version of the Loomis-Whitney inequality, known as the \emph{multilinear Kakeya} inequality, that was proved in the same paper \cite{bennett2006multilinear}. Next we explain this inequality.

\subsection{Multilinear Kakeya: a perturbed Loomis-Whitney inequality}\label{MKsec} If we take each $f_j$ to be a sum of finitely many characteristic functions of $1$-balls in the Loomis-Whitney inequality (\ref{gLW}), then each $f_j (x_1, \ldots, \hat{x}_j, \ldots, x_n)$ is a finite sum of characteristic functions $1_{T}$ where each $T$ is an infinite tube of radius $1$ (called unit tubes below) in the $x_j$ direction. From this viewpoint, the Loomis-Whitney inequality has a consequence:

\begin{thm}\label{LW2thm}
Fix the dimension $n$. Let $\Bbb{T}_j, 1\leq j \leq n,$ be a finite family of (infinite) unit tubes in $x_j$-direction. Then
\begin{equation}\label{cylinderLW}
    \int_{\Bbb{R}^n} \prod_{j=1}^n (\sum_{T \in \Bbb{T}_j} 1_T)^{\frac{1}{n-1}} \lesssim \prod_{j=1}^n |\Bbb{T}_j|^{\frac{1}{n-1}}.
\end{equation}
\end{thm}

We remark that Theorem \ref{LW2thm} is in fact equivalent to (\ref{gLW}), modulo a constant of loss. This can be seen by rescaling and an approximation argument. However, this version has a very intuitive geometric meaning that indicates we have a really sharp control on overlappings of cylinders from $n$ families in $n$ transverse directions. Moreover, if one tries to construct sharp examples for the multilinear restriction Theorem \ref{MRthm} from Knapp examples (i.e. single wave packets), like those in the discussion for the bilinear case, they quickly realize the alignment of the supports of wave packets (which are essentially tubes!) in those sharp examples will also give sharp examples for (\ref{cylinderLW}) after a rescaling. Therefore, it is very conceivable that close connections exist between the Loomis-Whitney inequality (\ref{cylinderLW}) and the multilinear restriction estimate (\ref{multilinearineq}).

\begin{figure}[htb]
		\centering
        \includegraphics[width=.2\textwidth]{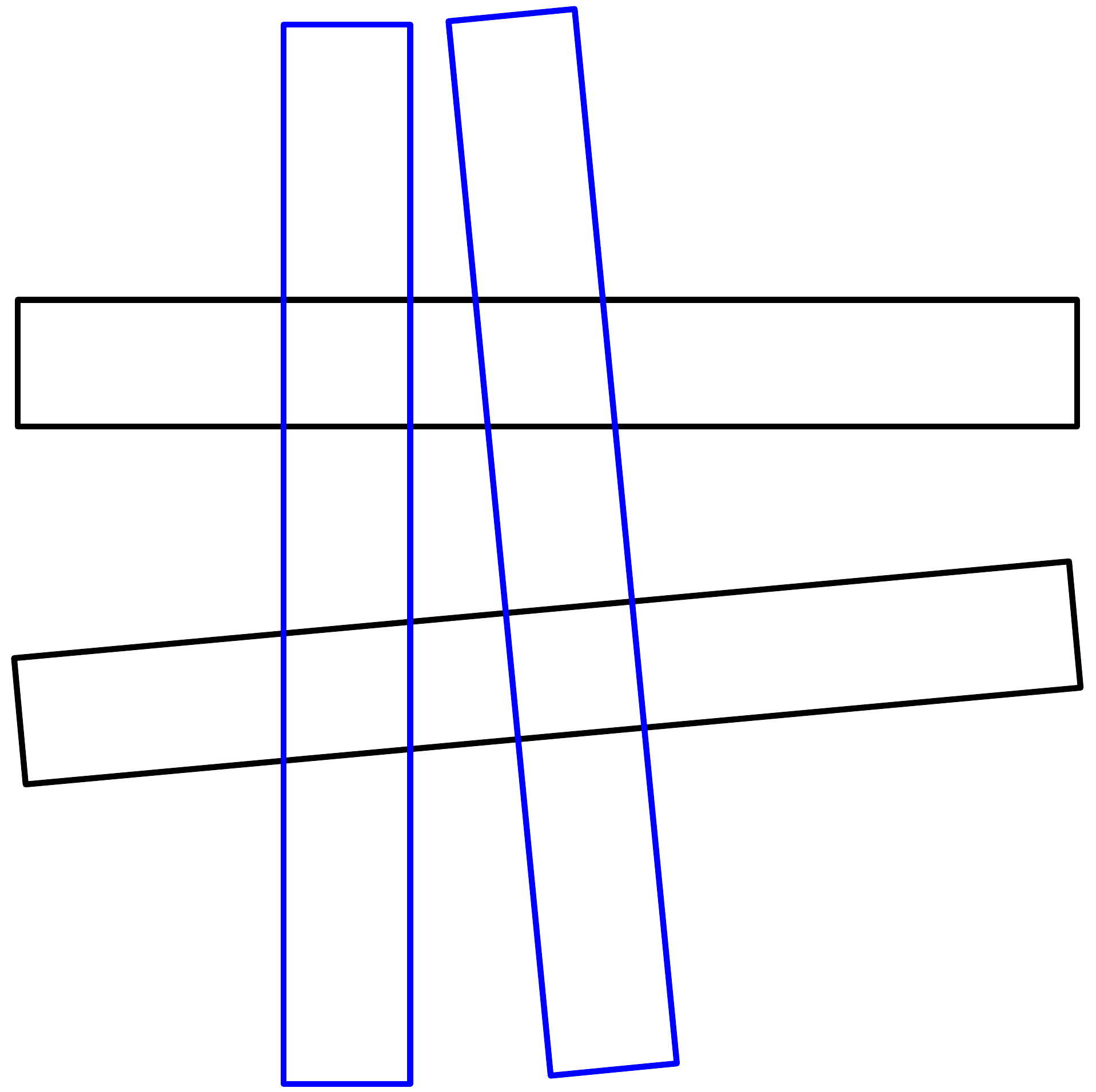}
        \centering\caption{The directions of wave packets in multilinear restriction may be tilted.}
        \label{MRpic}
\end{figure}

\begin{figure}[htb]
		\centering
        \includegraphics[width=.5\textwidth]{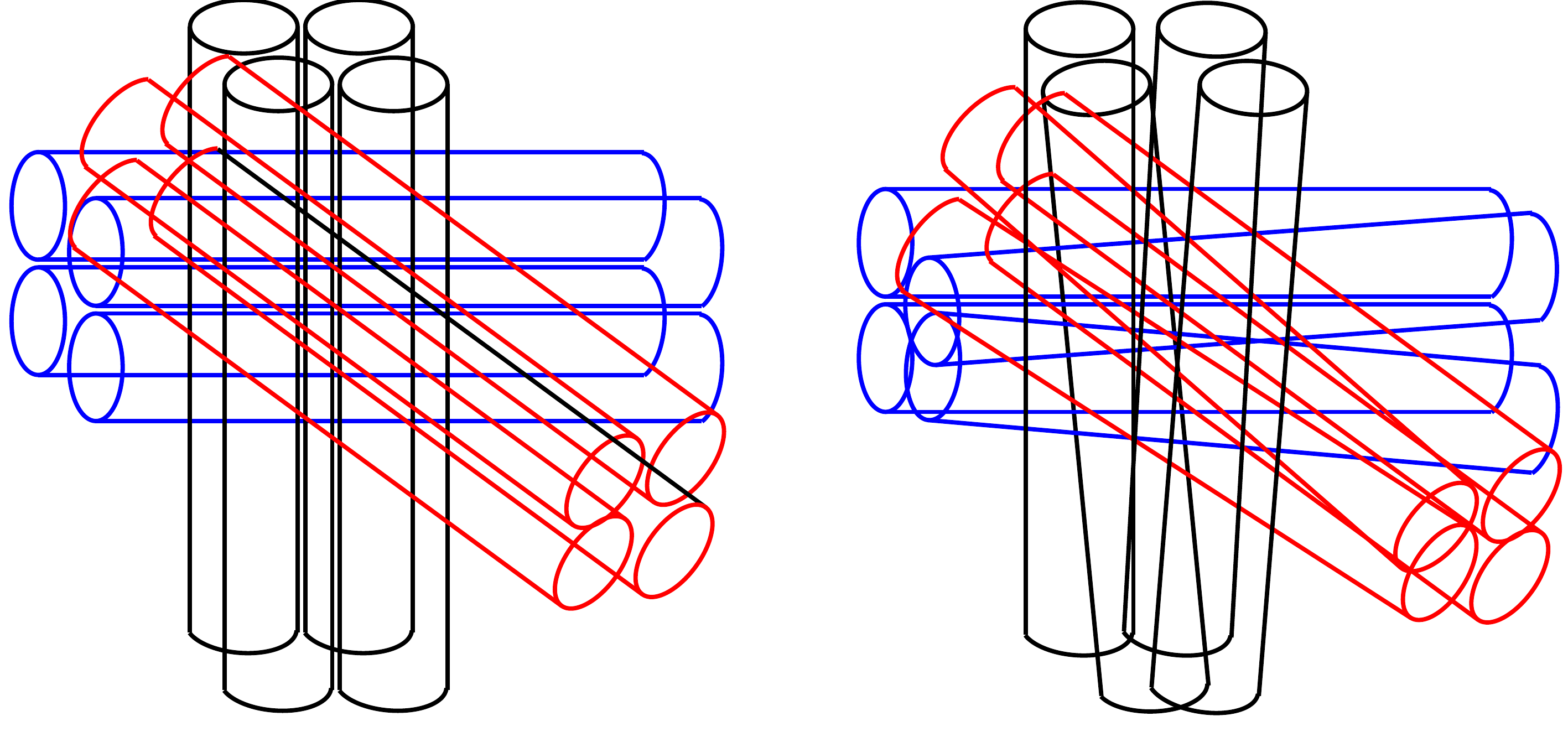}
        \centering\caption{Loomis-Whitney (left) and Multilinear Kakeya (right)}
        \label{LWvsMK}
\end{figure}

One then naturally wonders whether (\ref{cylinderLW}) can be used to prove (\ref{multilinearineq}). However, a glaring technical issue presents itself: For the multilinear restriction theorem there are other sharp examples where the main wave packets for each $f_j$ are essentially supported on disjoint $R^{\frac{1}{2}}-$tubes, but the directions of tubes slightly vary. See Figure \ref{MRpic}. In this setting, the amount of overlapping between different families of wave packet supports (i.e. supports of wave packets coming from $f_j$) should intuitively still be well-controlled but is not directly estimated by (\ref{cylinderLW}), since (\ref{cylinderLW}) only allow parallel tubes. This is the motivation of Bennett-Carbery-Tao when they developed a multilinear Kakeya estimate, which can be viewed as a much more robust version of (\ref{cylinderLW}) allowing ``perturbations'' of the directions of tubes within a family. See Figure \ref{LWvsMK} for an illustration.

\begin{thm}[Multilinear Kakeya \cite{bennett2006multilinear}]\label{MKthm}
Fix the dimension $n$, then there exists a small $c(n)>0$ such that: Let $\Bbb{T}_j, 1\leq j \leq n,$ be a finite family of (infinite) unit tubes such that each $T \in \Bbb{T}_j$ has an angle $\leq c(n)$ against the $x_j$-direction\footnote{We emphasize that different $T$ in the same $\Bbb{T}_j$ need not be parallel.}. Then inside every ball $B_R, R>1$,
\begin{equation}\label{MKineq}
    \int_{B_R} \prod_{j=1}^n (\sum_{T \in \Bbb{T}_j} 1_T)^{\frac{1}{n-1}} \lesssim_{\varepsilon} R^{\varepsilon} \prod_{j=1}^n |\Bbb{T}_j|^{\frac{1}{n-1}}, \forall \varepsilon > 0.
\end{equation}
\end{thm}

As before, restricting to a large ball $B_R$ and the loss $R^{\varepsilon}$ are usually completely acceptable in this kind of theorems.

\begin{rem}
It is standard to use ``Kakeya'' to name estimates that concern overlappings of tubes, after the important concept of Kakeya sets that we will not talk about much in this section. We will mention a bit about Kakeya conjectures in \S \ref{reallysharpsec}. For readers who are already familiar with Kakeya sets and Kakeya type estimates, it is worth noting that in Theorem \ref{MKthm}, tubes in a fixed family $\Bbb{T}_j$ are allowed to be parallel, which is forbidden in the classical linear Kakeya problem.
\end{rem}

The proof of Theorem \ref{MKthm} in \cite{bennett2006multilinear} is highly nontrivial and very unusual. Intuitively, they slightly enlarge the exponent by $\varepsilon$, then move all the tubes to center at the origin and verify that the left hand side will morally increase. Their approach is inspired by the aforementioned heat flow method in the proof of the Brascamp-Lieb inequality and a useful multiscale argument in Fourier analysis (known as ``induction on scales'', see also the next paragraph), but involves many new ideas including an ``extrapolation'' trick. Later Guth \cite{guth2015short} gave a shorter and conceptually simpler multiscale argument reproving Theorem \ref{MKthm}.

Once we have the more robust Multilinear Kakeya Theorem \ref{MKthm}, an iterative application of it will be enough to prove the Multilinear restriction Theorem \ref{MRthm}. To explain the proof, we notice that we are allowed to cut each $S_j$ into very tiny pieces and consider one piece at a time. Hence without loss of generality we may assume each $S_j$ has sufficiently small diameter (so that the normal vector from each $S_j$ will be a really tiny perturbation of a vector $\mathbf{v}_j$ and all $\mathbf{v}_j$ are quantitatively transverse). Then inside every ball $B_R$, wave packets from $f_j$ are really close to the $\mathbf{v}_j$ direction. Now one can use Theorem \ref{MKthm} to control their overlapping well. The rest is a typical ``induction on scales'' argument: By H\"{o}lder we have the desired conclusion for $R<100$ (say) when we take the implied constant sufficiently large. Now assuming we have the conclusion (``induction hypothesis'') on all $B_{R^{\frac{1}{2}}}$ in our $B_R$, we take a finitely overlapping covering of $B_R$ by balls $B_{R^{\frac{1}{2}}}$, one may prove the desired estimate (\ref{MRthm}) on $B_R$ by applying the induction hypothesis, the locally constant  and  local orthogonality properties of wave packets and finally applying (\ref{MKineq}).

\subsection{Perturbed Brascamp-Lieb inequalities}\label{PBLsec}

General Brascamp-Lieb inequalities (\ref{BL}), if true, also have corresponding perturbed versions. This was first proved by Bennett-Bez-Flock-Lee along with the proof of the stability Theorem \ref{BLstab}:

\begin{thm}[Perturbed Brascamp-Lieb inequality \cite{bennett2018stability}]\label{PBLthm}
Whenever (\ref{BL}) holds for some Brascamp-Lieb datum $(\mathbf{B}, \mathbf{p})$, there is a sufficiently small $c = c(\mathbf{B}, \mathbf{p}) > 0$ such that if we  take $\Bbb{T}_j$ to be $1$-neighborhoods of subspaces of $\Bbb{R}^n$ that form angles $<c$ against $\ker B_j$, then for every ball $B_R$, $R>1$, we have
\begin{equation}\label{PBLineq}
    \int_{B_R} \prod_{j=1}^m (\sum_{T \in \Bbb{T}_j} 1_T)^{p_j} \lesssim_{\varepsilon} R^{\varepsilon} \prod_{j=1}^m |\Bbb{T}_j|^{p_j}, \forall \varepsilon > 0.
\end{equation}
\end{thm}

The importance of (\ref{PBLineq}) is natural in Fourier restriction theory for objects other than the paraboloid or the sphere, such as the moment curve. As one can imagine, for a general submanifold, small caps on it are possibly better approximated by boxes that are not necessarily ``slabs''. For example each $\delta$-segment on the moment curve $\mathcal{M}_3$ can be put in a box of dimensions $\sim\delta \times \delta^2 \times \delta^3$. This suggests that the locally constant regions as well as the wave packet decomposition can all be in different shapes from tubes. Indeed, Theorem \ref{PBLthm} plays a crucial role in Bourgain-Demeter-Guth's proof \cite{bourgain2016proof} of the sharp solution count of system (\ref{Vleeqndeg4}). Starting in the next section, we will talk more about their approach as one of several examples of applications of Theorems \ref{MKthm}, \ref{MRthm} and \ref{PBLthm} in Fourier restriction type problems.

\subsection{Additional historical remarks}\label{AHR3} The Stein-Tomas theorem is closely connected to the very important Strichartz estimates in PDE \cite{strichartz1977restrictions, keel1998endpoint}. See e.g. Lecture 6 of \cite{tao2004recent} for Strichartz estimates, their applications to nonlinear PDE and more. These estimates were recently generalized to ones for systems of orthonormal functions by Frank and Sabin \cite{frank2017restriction} that are important in the context of many-body quantum physics. See also Frank-Lieb-Sabin \cite{frank2016maximizers} for discussions of maximizers of the Stein-Tomas inequality.

Bilinear methods originated from classical $L^4$ estimates (\cite{fefferman1973note, carleson1972oscillatory, cordoba1977kakeya, carbery1983boundedness, mockenhaupt1993note}). Since $4$ is an even number, it is very natural to use Plancherel to study $L^4$ estimates. In such an approach to the two dimensional Fourier extension problem, one then sees that curvature conditions are translated to transversality conditions for far-apart pieces on the curve, which then motivates the bilinear approach (see e.g. Section 2.2 of \cite{tao2004recent} for a more detailed discussion). For more recent bilinear methods in Fourier restriction type problems (with the exponent $<4$), see \cite{bourgain1995estimates, tao1998bilinear, tao2000bilinear, tao2000bilinear2, wolff2001sharp, tao2003sharp}.

For the origin and motivations of $\varepsilon$-removal lemmas, see Theorem 2.11 in \cite{tao2004recent}, the paragraphs around it, and e.g. the references   therein (like \cite{bourgain1995estimates, tao1998weak, tao1999bochner, tao2000bilinear}).

For some earlier applications of induction on scales in Fourier analysis, see \cite{bourga1991besicovitch, wolff2001sharp, tao2003sharp}.

\section{Applications to decoupling theory}

One and a half decades have passed after the multilinear Kakeya Theorem \ref{MKthm} was originally discovered. Today that theorem and the generalized Theorem \ref{PBLthm} are ubiquitous in the study of Fourier restriction problems. In this section we sketch a very typical and influential example: Applications of Theorems \ref{MKthm} and \ref{PBLthm} in proofs of a class of important tools in Fourier restriction theory known as \emph{decoupling}. A good example of references on decoupling is Demeter's book \cite{demeter2020fourier}.

\subsection{Decoupling estimates} Decoupling estimates are a kind of useful estimates that were first introduced by Wolff \cite{wolff2000local}, where he used these estimates to make progress on the local smoothing Conjecture \ref{localsmoothingconj}. We begin by stating the important decoupling theorem for $P^{n-1}$ of Bourgain-Demeter \cite{bourgain2015proof}. Take an arbitrarily large parameter $R>1$ and partition  $P^{n-1}$ into disjoint $R^{-\frac{1}{2}}$-caps $\theta$ as before. We are interested in ways to estimate a function $g$ on a ball $B_R$ with the Fourier support of $g$ in $P^{n-1}$.  Define $\hat{g}_{\theta}$ to be the restriction of $\hat{g}$ on $\theta$ as before, then we have the following: 

\begin{thm}[Decoupling, Bourgain-Demeter \cite{bourgain2015proof}]\label{BDthm}
For $2 \leq p \leq \frac{2(n+1)}{n-1}$,
\begin{equation}\label{BDthmineq}
    \|g\|_{L^p (B_R)} \lesssim_{\varepsilon} R^{\varepsilon} (\sum_{\theta} \|g_{\theta}\|_{L^p (w_{B_R})}^2)^{\frac{1}{2}}, \forall \varepsilon>0.
\end{equation}
\end{thm}

As a technical setup, $w_{B_R}$ is a weight that is $1$ on $B_R$ and decays sufficiently rapidly away from $B_R$. Morally speaking, (\ref{BDthmineq}) is almost as strong as the following square-root cancellation type inequality
\begin{equation}\label{BDthmfakeineq}
    \|g\|_{L^p (B_R)} ``\lesssim" (\sum_{\theta} \|g_{\theta}\|_{L^p (B_R)}^2)^{\frac{1}{2}}.
\end{equation}
After Theorem \ref{BDthm}, it is a good intuition to pretend (\ref{BDthmfakeineq}) holds in applications. And when arguing rigorously one instead uses the inequality (\ref{BDthmineq}) which has a weight $w_{B_R}$ with a tail for purely technical reasons. Since we will mainly be talking about intuitions, we will often not distinguish between $L^p (B_R)$ and $L^p (w_{B_R})$ as the tail of the weight only comes as a technical consequence of the rigorous form of the uncertainty principle and ignoring it will not affect the intuition in proofs. In view of this, we will often not distinguish between, e.g., (\ref{BDthmineq}) and (\ref{BDthmfakeineq}).

Arguably, the decoupling Theorem \ref{BDthm} is perhaps the single most influential theorem in the past decade in Fourier restriction theory. The theorem has important direct consequences such as improved discrete restriction and the optimal Strichartz estimates for the Schr\"{o}dinger equation on the torus (see Section 2 in \cite{bourgain2015proof} or Section 13.3 in \cite{demeter2020fourier}). For contexts of these problems see for example \cite{bourgain1993fourier} or Chapter 1 in \cite{bourgain2007mathematical}. Moreover, Theorem \ref{BDthm} played an important role in the proofs of breakthroughs such as new results on the Fourier restriction/extension Conjectures \ref{restrconj} and \ref{extconj} by Guth \cite{guth2018restriction} and subsequent works, a complete proof of the Main Conjecture in Vinogradov's Mean Value Theorem \cite{bourgain2016proof} that in particular gives optimal solution count of the system (\ref{Vleeqndeg4}), and all the aforementioned recent progress on Falconer's distance Conjecture \ref{Falconerconj}. Its proof also inspired new methods that, for example, lead to the proof of the local smoothing Conjecture \ref{localsmoothingconj} in dimension $2+1$ by Guth, Wang and the author \cite{guth2020sharp}.

\subsection{The role of multilinear Kakeya in the proof} Bourgain and Demeter's proof of the decoupling Theorem \ref{BDthm} has two main components. One is the application of the multilinear Kakeya Theorem \ref{MKthm} and the other is a series of carefully designed arguments that induct on scales. Since the second component is very sophisticated, we will not be able to even give a description of the proof on an intuitive level. Instead, in this subsection we try to intuitively describe the key role played by multilinear Kakeya in Bourgain and Demeter's proof.

When trying to prove an inequality like (\ref{BDthmfakeineq}), it is very helpful to define the following quantity and track how it changes with $r$ and $\sigma$ increasing dyadically from $1$:
\begin{equation}\label{Dpq}
    D_{p, q} (r, \sigma) = \sum_{B_r \subset B_R} (\sum_{\tau: \sigma^{-1}-\text{cap}} \|g_{\tau}\|_{L^q (B_r)}^2)^{\frac{p}{2}}
\end{equation}
where $p, q \geq 2$ 
and the collection of $B_r \subset B_R$ always means a finitely overlapping covering. Indeed, by the uncertainty principle, when $r = \sigma =1$, (\ref{Dpq}) is morally the left hand side of (\ref{BDthmfakeineq}) raised to the $p$-th power, while when $q=p$, $r= R$ and $\sigma = R^{\frac{1}{2}}$, (\ref{Dpq}) is morally the right hand side of (\ref{BDthmfakeineq}) raised to the $p$-th power. The strategy of Bourgain-Demeter's proof focuses on this quantity $D_{p, q}$ and related quantities. The parameters $r$ and $\sigma$ are made to grow gradually while careful analysis is made so that we almost do not lose any factor along the way. Note that in the whole game we are only allowed to lose an $R^{\varepsilon}$.

Multilinear Kakeya is the crucial ingredient in the step when $r$ is increased (known as ``ball inflation'' in the literature). To explain this, let us first make some preliminary observations. 

\begin{figure}[htb]
		\centering
        \includegraphics[width=.2\textwidth]{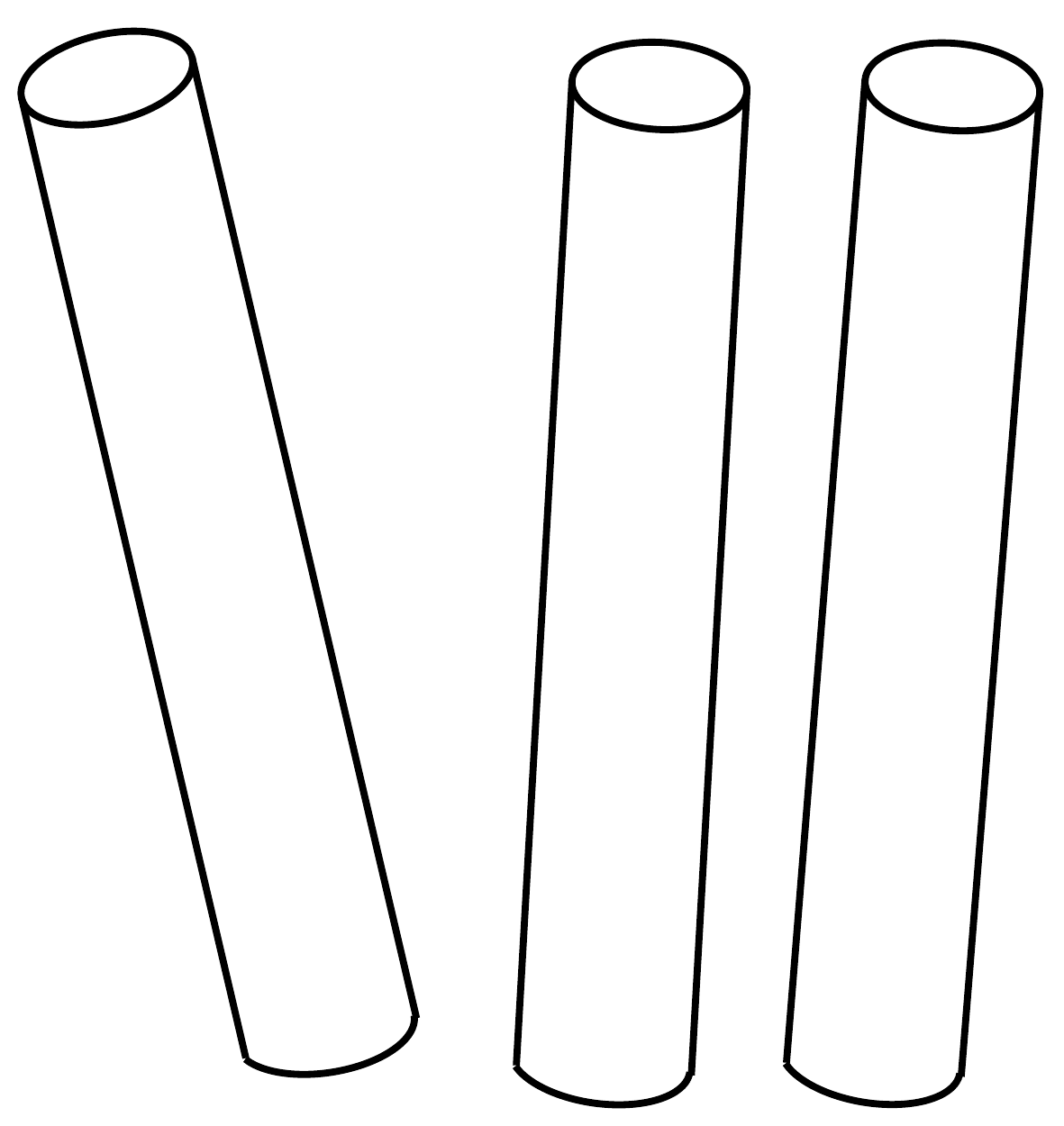}
        \centering\caption{The disjoint supports case.}
        \label{Disjoint}
\end{figure}

First, we see that Theorem \ref{BDthm} is trivially true if all but one $g_{\theta}$ is zero (and hence $g = g_{\theta}$). More generally, Theorem \ref{BDthm} is morally straightforward if the essential support\footnote{Recall from the locally constant property that each $|g_{\theta}|$ is like a constant on tubes of length $R$ and thickness $R^{\frac{1}{2}}$. So in general one can imagine each $g_{\theta}$ is essentially supported on the union of a few parallel, finitely-overlapping tubes.} of all $g_{\theta}$ inside $B_R$  are disjoint (Figure \ref{Disjoint}), because in this case at each individual point in $B_R$, $g$ will essentially be equal to one $g_{\theta}$.

\begin{figure}[htb]
		\centering
        \includegraphics[width=.3\textwidth]{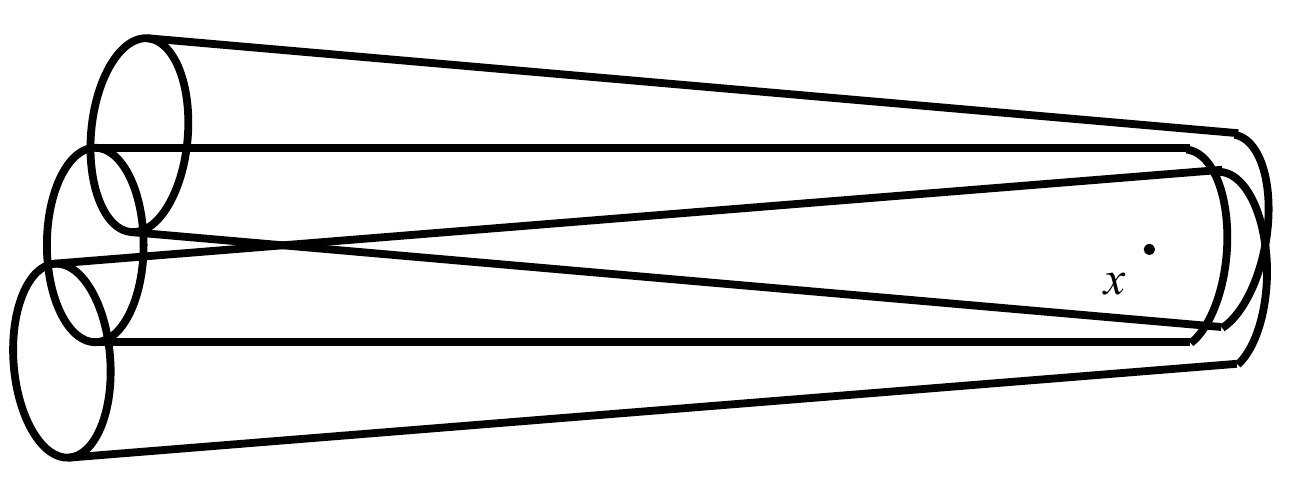}
        \centering\caption{The narrow case.}
        \label{Narrow}
\end{figure}

An observation by Bourgain-Guth \cite{bourgain2011bounds} in a different context allows us to deal with more cases by induction on dimension, such as the following ``narrow'' case: at each point $x$ in $B_R$, the main contribution to $g(x)$ comes from $g_{\theta} (x)$'s whose wave packets are all parallel to one hyperplane (Figure \ref{Narrow}). By Bourgain-Guth's method, one can induct on dimension to use Theorem \ref{BDthm} in dimension $n-1$ to prove Theorem \ref{BDthm}.\footnote{Strictly speaking, one also induct on the scale $R$. But we will not go to full details.}

\begin{figure}[htb]
		\centering
        \includegraphics[width=.3\textwidth]{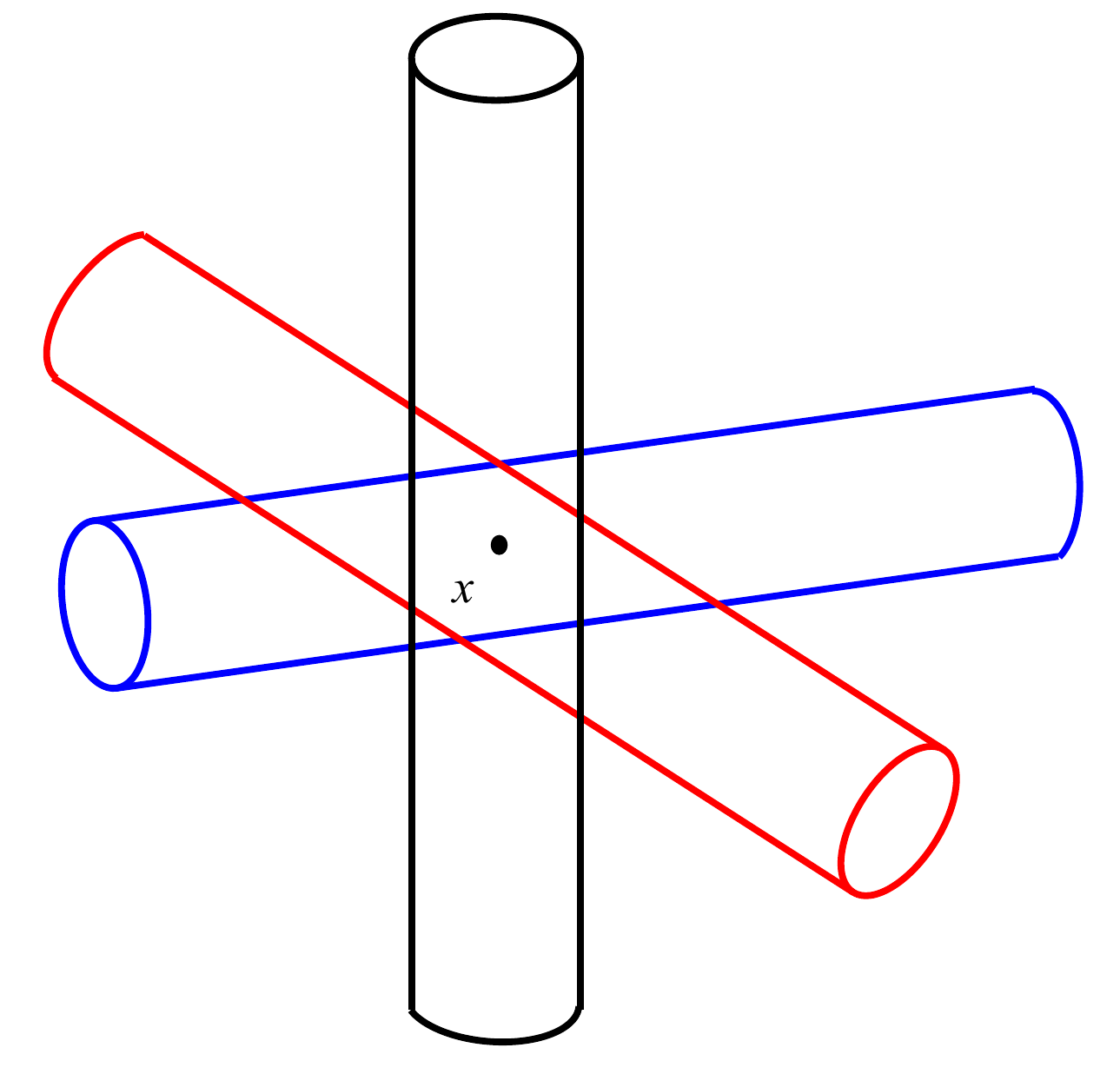}
        \centering\caption{The broad case.}
        \label{Broad}
\end{figure}

We do not plan to give the full Bourgain-Guth method here, but at this point will just say that in light of this reduction, we can further assume for the most interesting case of Theorem \ref{BDthm} that, at a typical point $x$ in $B_R$, we can find $n$ transverse directions such that for every such direction, the $g_{\theta}$ whose wave packets are ``close'' to that direction together make significant contribution to $g(x)$ (Figure \ref{Broad}). When this happens, we can bound $D_{p, q} (r, r)$ by $D_{p, q} (r^2, r)$ with almost no loss by the multilinear Kakeya Theorem \ref{MKthm} for $q = \frac{n}{n-1} p$. For readers not working in decoupling but wish to learn Bourgain-Demeter's proof, it is a very good exercise to work out the details of both the statement and the proof of this ``ball inflation'' argument (which are beyond our scope here). As a hint we just remark that by the locally constant property, each $g_{\tau}$ will behave like a  constant on tubes of length $r^2$ and thickness $r$ and hence by definition (\ref{Dpq}) this bound is really a geometric combinatorial inequality. Multilinear Kakeya will fit perfectly into the picture with our   ``transverse contributions from wave packets'' assumption in the beginning of this paragraph.

\begin{rem}
In completely rigorous terms, to prove Theorem \ref{BDthm} one needs to first formulate a multilinear version and consider the multilinear version of $D_{p, q}$:
\begin{equation}\label{Mpq}
    M_{p, q} (r, \sigma) = \text{Avg}_{B_r \subset B_R} \prod_{j=1}^n(\sum_{\tau \subset \Omega_j: \sigma^{-1}-\text{cap}} \|g_{\tau}\|_{L_{\text{Avg}}^q (B_r)}^2)^{\frac{p}{2n}}
\end{equation}
where Avg stands for average and $\Omega_1, \Omega_2, \ldots, \Omega_n$ are regions on $P^{n-1}$ such that they have quantitatively transverse normal directions. Multilinear Kakeya is then used to prove the ball inflation lemma
\begin{equation}
    M_{p, \frac{n-1}{n} p} (r, r) \lesssim_{\varepsilon} r^{\varepsilon} M_{p, \frac{n-1}{n} p} (r^2, r).
\end{equation}

The Bourgain-Guth argument that we introduced above is used to show that the multilinear version implies the original version of Theorem \ref{BDthm}.
\end{rem}

\begin{rem}\label{BGoriginal}
Originally, Bourgain-Guth developed their method in \cite{bourgain2011bounds} to go from multilinear Kakeya to new honest Fourier restriction estimates. This is another big application of multilinear Kakeya in its early stage. It transpired that the nature of the Fourier restriction conjectures is slightly different from that of decoupling in the following sense: The most difficult scenario in the Fourier restriction conjectures does not appear to be the case where at every point the main contributions are ``transverse''. For example in dimension $3$, a pretty difficult scenario is that we have ``planar'' contribution of wave packets around each point, i.e. the main contribution to every point comes from wave packets near one plane (that may rotate depending on the point). For this reason, the Fourier restriction conjectures \ref{restrconj} and \ref{extconj} are still widely open and it is generally believed that the strength of multilinear Kakeya/restriction alone is not enough for full resolutions.
\end{rem}

\subsection{Applications in the Main Conjecture Vinogradov's Mean Value Theorem} The Main Conjecture in Vinogradov's Mean Value Theorem predicts a sharp bound of the count of integer solutions to (\ref{Vleeqndeg4}) and similar systems in other degrees. See \cite{wooley2014translation} for an excellent introduction. When the degree is $\leq 2$, the conjecture is elementary and easy. Wooley \cite{wooley2016cubic} proved the cubic case of the Main Conjecture and Bourgain-Demeter-Guth \cite{bourgain2016proof} prove the Main Conjecture for all degrees.

Bourgain-Demeter-Guth's result is obtained via a decoupling result in the form similar to (\ref{BDthmineq}) or (\ref{BDthmfakeineq}) (but with a different exponent) concerning partition of the moment curve $\mathcal{M}_n$ into segments $J$ of length $\delta$. If $\text{supp } \hat{g} \subset \mathcal{M}_n$ then we let $\hat{g}_J$ to be the restriction of $\hat{g}$ on $J$. Bourgain-Demeter-Guth deduced the Main Conjecture from
\begin{thm}[Decoupling for the moment curve \cite{bourgain2016proof}]\label{BDGthm}
For $2 \leq p \leq n(n+1)$,
\begin{equation}\label{BDGthmineq}
    \|g\|_{L^p (B_{\delta^{-n}})} \lesssim_{\varepsilon} R^{\varepsilon} (\sum_{J} \|g_{J}\|_{L^p (w_{B_{\delta^{-n}}})}^2)^{\frac{1}{2}}, \forall \varepsilon>0.
\end{equation}
\end{thm}

In their proof, the perturbed Brascamp-Lieb inequalities (cf. Theorem \ref{PBLthm})  played a key role. To see this, we first mention that one can again use Bourgain-Guth style arguments to ``multilinearize'' the problem, then rely on ``ball-inflation'' techniques. What is different here is that they are able to utilize a whole lot of different ball inflation inequalities (used to inflate the ball and increase the radius from $r^m$ to $r^{m+1}$ with $m = 1, 2, \ldots, n-1$). As before, these inequalities eventually come from the locally constant properties of functions such as $g_J$, as we briefly explain next.

Elementary differential geometry shows that each $J$ can be put in a rectangular box of dimensions $\sim \delta \times \delta^2 \times \cdots \times \delta^n$. Hence by the uncertainty principle, each $|g_J|$ is like a constant on boxes of dimensions $\sim \delta^{-1} \times \delta^{-2} \times \cdots \times \delta^{-n}$. Why is this relevant? Let us take any $m$ between $1$ and $n-1$. For each point $x$ in $B_{\delta^{-n}}$, let $F_{m, J} (x)$ be the $L^q$ integral of $|g_J|$ on a ball of radius $\delta^{-m}$ around $x$ which is one of the quantities coming up in the analogue of (\ref{Dpq}). Because of the locally constant property, we see
\begin{prin}\label{multiuncertaintyprin}
$F_{m, J}$ behaves like a constant on every translation of a box $\Box_J$ of dimension $\delta^{-m} \times \cdots \times \delta^{-m} \times \delta^{-m-1} \times \cdots \times \delta^{-m-1}$, where $m$ of the edges have length $\delta^{-m}$, and the direction of $\Box_J$ changes depending on $J$.
\end{prin}
To effectively control the transverse intersections of boxes with multiple longest edges, perturbed Brascamp-Lieb inequalities instead of multilinear Kakeya needs to be used. Bourgain-Demeter-Guth showed that for every $m$, as long as one takes a sufficiently large number of boxes that correspond to segments in ``generic positions'', they indeed have good perturbed Brascamp-Lieb inequalities needed (stated and proved in Section 6 in \cite{bourgain2016proof}). This can be proved using the ``non-degenerate'' way in which $\mathcal{M}_n$ is ``curved''. It is then a combination of this geometric fact and an involved induction on scales machinery that enabled them to solve the full Main Conjecture.

\begin{rem}
Wooley's proof \cite{wooley2016cubic} in the degree $3$ case and his later proof of the Main Conjecture in the full generality \cite{wooley2019nested} use a different technique known as efficient congruencing. This method and the method of decoupling have a lot of components analogous to each other but they also have different aspects.
\end{rem}

\subsection{Additional historical remarks}\label{AHR} 

The idea of decoupling can be traced back to Wolff's paper \cite{wolff2000local}. Early developments involve extensions to higher dimensions by \L{}aba–Wolff \cite{laba2002local} and refinements by Garrig\'{o}s–Seeger \cite{garrigos2009plate, garrigos2010mixed} and Garrig\'{o}s–Schlag–Seeger \cite{garrigos2008improvements}. Bourgain also had a partial result \cite{bourgain2013moment} towards Theorem \ref{BDthm}. The solutions to Carleson's pointwise convergence problem (introduced in \cite{carleson1980some}) for free Schr\"{o}dinger solutions in higher dimensions \cite{du2017sharp, du2019sharp} are another two important recent results whose proofs utilize both decoupling and multilinear Kakeya as key tools.

\section{Decoupling for TDI systems} In the latest development of Fourier restriction theory, decoupling and perturbed Brascamp-Lieb inequalities have continued to play fundamental parts. In this section we talk about a problem that generalizes the Main Conjecture in Vinogradov's Mean Value Theorem and has seen a lot of recent progress -- solution counting for TDI (translation-dilation-invariant) systems, defined in (\ref{TDIsystemeq}) below. Recently several very satisfactory results were obtained and it turns out that one of the most recent developments on this problem utilizes the perturbed regularized Brascamp-Lieb inequality (or called perturbed scale-dependent Brascamp-Lieb inequality), which is a further generalization of Theorem \ref{PBLthm}.

\subsection{Solution counting for TDI systems} Recall that Vinogradov's idea of bounding the number of solutions to (\ref{Circleeqndeg4}) when the number of variables is sufficiently large is to add auxiliary equations to arrive at system (\ref{Vleeqndeg4}). We mentioned in \S\ref{FRTprob} that system  (\ref{Vleeqndeg4}) behaves better because it is connected to the moment curve $\mathcal{M}_4$ with better ``curvature'' conditions. Now we make an observation that justifies one aspect of this ``goodness''. 

Note that $\mathcal{M}_4$ (or general $\mathcal{M}_n$) has a good property that any segment on it can be rescaled under an affine transform to be the whole curve. For details of this property, we refer the readers to e.g. the one dimensional analogue of the discussion in Section 7 in \cite{bourgain2016mean}). Intuitively speaking, this means that the way $\mathcal{M}_4$ is curved is  ``very similar'' in all places and has ``self-similarity''. This symmetry property is known as the \emph{translation-dilation invariant (TDI) property} of the curve. A great number of submanifolds in $\Bbb{R}^n$ have this property, including the paraboloid $P^{n-1}$ that we have been focusing on. The transform that takes one small piece of $P^{n-1}$ to the whole $P^{n-1}$ is known as a \emph{parabolic rescaling} and plays a fundamental role in every approach to the Fourier extension Conjecture \ref{extconj} that uses induction on scales. We will see more submanifolds of $\Bbb{R}^n$ with the TDI property below. They will be called \emph{TDI manifolds}.

On the circle method side, it turns out that there is a concept of \emph{TDI systems} that are related to TDI manifolds: A TDI system has a related TDI manifold, whose geometric properties are strongly related to the solution counting of the system. TDI systems look like:
\begin{equation}\label{TDIsystemeq}
    \sum_{j=1}^s
(\mathbf{F} (\mathbf{x}_j) - \mathbf{F} (\mathbf{y}_j)) = \mathbf{0} \end{equation}
where each component $F_h$ of $\mathbf{F}$ is a homogeneous form and the system is to be solved for $|x_{jk}| , |y_{jk}| \leq N$ and $x_{jk} , y_{jk} \in \Bbb{Z}$ with $N \to \infty$. Additionally, the TDI property translates to the following property of (\ref{TDIsystemeq}): for each $F_h$, each of its partial derivatives is a linear combination of some components of $\mathbf{F}$.

\begin{eg}
The following system of five equations
\begin{equation}
\begin{cases}
      x_{1, 1}+ x_{2, 1} +x_{3, 1} = y_{1, 1} + y_{2, 1} + y_{3, 1},\\
      x_{1, 2}+ x_{2, 2} +x_{3, 2} = y_{1, 2} + y_{2, 2} + y_{3, 2},\\
      x_{1, 1}^2+ x_{2, 1}^2 +x_{3, 1}^2 = y_{1, 1}^2 + y_{2, 1}^2 + y_{3, 1}^2,\\
      x_{1, 2}^2+ x_{2, 2}^2 +x_{3, 2}^2 = y_{1, 2}^2 + y_{2, 2}^2 + y_{3, 2}^2,\\
      x_{1, 1} x_{1, 2}+ x_{2, 1} x_{2, 2} + x_{3, 1} x_{3, 2} = y_{1, 1} y_{1, 2}+ y_{2, 1} y_{2, 2} + y_{3, 1} y_{3, 2}\\
\end{cases}
\end{equation}
is a TDI system.
\end{eg}

Given a TDI system we can ask the following:

\begin{ques}[Sharp solution counting for TDI systems]\label{solucountques}
Given a TDI system in the shape of (\ref{TDIsystemeq}) to be solved for $|x_{jk}| , |y_{jk}| \leq N$ and $x_{jk} , y_{jk} \in \Bbb{Z}$. What is the largest $\beta$ such that the number of solutions to the system is $O_{\varepsilon} (N^{\beta+\varepsilon}), \forall \varepsilon > 0$?
\end{ques}

Because of the symmetry of TDI manifolds, additional tools are available for bounding the number of solutions to TDI systems and  there is now serious hope to solve Question \ref{solucountques} in great generality. These systems therefore attracted a large amount of recent interest.

Some general and sharp solution-counting results are already achieved by the method of decoupling. Indeed, as with the Vinogradov system, TDI systems can be studied via decoupling for the corresponding TDI manifolds. Next we explain this approach in more details.

\begin{rem}
We only focus on the method of decoupling in this article, but it should be pointed out that the method of efficient congruencing is also quite influential and powerful and many connections made between people's work with both methods have in turn benefited each approach a lot.
\end{rem}

\subsection{Sharp decoupling for TDI manifolds} We define a TDI manifold to be a submanifold of $\Bbb{R}^{d+n}$ that is a graph defined by \[\xi_{d+j} = F_{j} (\xi_1, \ldots, \xi_d), 1\leq j \leq n\] where each $F_j$ is a homogeneous polynomial whose all partial derivatives are linear combinations of $F_{1}, \ldots, F_{n}$.\footnote{One can use more intrinsic definitions but the definition here already includes all interesting examples and has the advantage of being concrete.} 
One can check such manifolds have very good symmetries (translation-dilation invariance) as mentioned above.

Let $S$ be the part of a TDI manifold satisfying $|\xi_i| \leq 1, \forall 1 \leq i \leq d$. For a parameter $\delta < 1$, divide $S$ into \emph{$\delta$-caps} $\theta$ defined by \[\xi \in S: M_i \delta \leq \xi_i < (M_i+1) \delta, \forall 1 \leq i \leq d\] where each $M_i$ runs through all integers.\footnote{Note there are only finitely many (about $\delta^{-d}$) non-empty caps. We will drop all empty ones.} In applications to solution counting of TDI systems, it is often useful to look for \emph{sharp decoupling} inequalities for $S$ of the following form, which is an important Fourier restriction type problem about $S$ in its own right.

When $g$ is a function on $\Bbb{R}^{d+n}$ whose Fourier support is in $S$, we use the notation that $\hat{g}_{\theta}$ is the restriction of $\hat{g}$ on $\theta$. We are then interested in the sharp decoupling inequality:
\begin{equation}\label{decgeneral}
    \|g\|_{L^p (B_{\delta^{-k}})} \lesssim_{\varepsilon} R^{\gamma + \varepsilon} (\sum_{\theta} \|g_{\theta}\|_{L^p (w_{B_{\delta^{-k}}})}^q)^{\frac{1}{q}}, \forall \varepsilon>0,
\end{equation}
where $\gamma = \gamma (p, q, S) \geq 0$, the \emph{degree} $k$ is defined as the highest degree of all $F_j$, and  being sharp simply means that the inequality would fail if $\gamma$ is replaced by a smaller non-negative real number. We will call (\ref{decgeneral}) an \emph{$l^q L^p$-decoupling} inequality.

Very often, sharp solution counting for TDI systems follows from some sharp decoupling theorem. Recall one example that Bourgain, Demeter and Guth \cite{bourgain2016proof} deduced the Main Conjecture in Vinogradov's Mean Value Theorem for arbitrary degree using the sharp $l^2 L^p$-decoupling Theorem \ref{BDGthm}. The sharp decoupling Theorem \ref{BDthm} is also $l^2 L^p$-decoupling. However, it is known that in order to get sharp solution counting for general TDI systems, one is sometimes forced to focus on proving sharp non-$l^2L^p$-decoupling inequalities due to an ``intermediate dimension'' issue. In other words, one has to consider the case $q>2$ in some cases.

Recently, there has been a lot of success in proving more sharp decoupling inequalities and using them to verify sharp solution counting for corresponding TDI systems. In the next two sections, we mention two such examples -- sharp decoupling for Parsell-Vinogradov manifolds and for arbitrary TDI manifolds of degree $2$. As we will see, in both examples, new Brascamp-Lieb-related ingredients are developed and become one type of central components in the proof.

\subsection{Proving finiteness of the Brascamp-Lieb constant in a systematic way}\label{systemBL} As in the cases of paraboloids and moment curves, for a general TDI manifold $S$, by a Taylor series expansion argument, we see each $\delta$-cap can be put in a rectangular box, each of whose edge length is $\sim \delta, \sim \delta^2, \ldots,$ or  $\sim\delta^k$. The locally constant property suggests that each $|g_{\theta}|$ is like a constant on the corresponding dual box, which has dimensions all roughly being positive integer powers of $\delta^{-1}$ (up to $\delta^{-k}$). In particular, as we did in the analysis for the moment curve, for any given $1 \leq j \leq k-1$, $|g_{\theta}|$ is essentially a constant on all translations of a particular box, each of whose edge lengths is either $\sim \delta^{-1}$ or $\sim \delta^{-j-1}$. The direction of the box depend on $\theta$. As we saw from the two previous examples (Theorems \ref{BDthm} and \ref{BDGthm}), proving different such boxes (in more rigorous terms, thickenings of those boxes by $\delta^{-j+1}$ times) cannot overlap a lot is one of the keys in proofs of decoupling inequalities for $S$. To prove such a statement, one can imagine perturbed Brascamp-Lieb inequalities to help.

Indeed, one common feature shared by the decoupling for the paraboloid (Theorem\ref{BDthm}) and the decoupling for the moment curve (Theorem \ref{BDGthm}) is that there is some useful Brascamp-Lieb inequality for a sufficiently large number of ``generic'' boxes.\footnote{To avoid further technicality, we choose not to explain the precise meaning of ``genericity'' here but only mention that the seemingly-bad non-generic case can always be taken care of by a Bourgain-Guth style argument.} This ``non-perturbed'' version then implies a useful perturbed Brascamp-Lieb inequality (by Theorem \ref{PBLthm}), which is then a crucial ingredient\footnote{We remark that other significant difficulties, mainly in the design of the machine that does induction on scales sharply, also need to be overcome. But we do not focus on them in this article.} of the proof of sharp decoupling. This feature continues to hold for a great number of other TDI manifolds. We mention one example here: a higher dimensional \emph{Parsell-Vinogradov manifold} of dimension $d$ and degree $k$ is defined as a joint graph of all non-constant monomials in $d$ variables of degree $\leq k$: \[(\xi_1, \xi_2, \ldots, \xi_d, \xi_1^2, \xi_1 \xi_2, \ldots, \xi_1 \xi_d, \xi_2^2, \ldots, \xi_d^2, \xi_1^3, \xi_1^2 \xi_2, \ldots, \xi_d^3, \ldots, \xi_1^k, \xi_1^{k-1} \xi_2, \ldots, \xi_d^k).\]

These (TDI) manifolds are higher dimensional generalizations of the moment curves and their sharp decoupling is now known ($d=1$ by \cite{bourgain2016proof}, $d=2, k=3$ by \cite{bourgain2017sharp} and in the full generality by \cite{guo2019integer}). For all of these manifolds, it turns out that for sufficiently many generic caps, one has useful (non-perturbed and perturbed) Brascamp-Lieb inequalities. More precisely, this means that when one looks at boxes of various scales, each coming from a different cap (similar to those in Principle \ref{multiuncertaintyprin}), as long as there are sufficiently many caps in generic positions, the translations of such boxes (together with some choices of exponents) satisfy the Brascamp-Lieb inequality, in the sense explained in \S \ref{MKsec} and Theorem \ref{PBLthm}. 

Assuming there are $M$ caps, one sees from differential geometry that having the Brascamp-Lieb inequality above precisely means that for every $1 \leq h< k$ there exist some $p>0$ such that when we denote $H_{h, \theta}$ to be the $h$-th tangent space (in geometric terms, these are just images of the $h$-th jet spaces) at the center $c_{\theta}$ of a cap $\theta$ and $B_{h, \theta}$ to be the orthogonal projection to $H_{h, \theta}$, we have a finite Brascamp-Lieb constant for the Brascamp-Lieb datum $(B_{1, \theta}, \ldots,  B_{M, \theta}, p, \ldots, p)$. We remark that one does not lose much by restricting their attention to the current setting where all exponents are the same (all being $p$) by noticing the natural symmetry in the concept of genericity.
 
This can be made very concrete with the help of Theorem \ref{BCCTfinite}. Denote the total dimension of the ambient space to be $d+n = {{d+k}\choose d} -1$ and the dimension of each $h$-th tangent space by $d_h$ (note that because of the translation invariance, all $h$-th tangent spaces have the same dimension). We need no more and no less than the following:

(i) $n+d = d_h Mp$

(ii) $\dim V \leq p \sum_j \dim B_{j, \theta} V$, $\forall$ subspace $V \subset \Bbb{R}^{d+n}$.

Eliminating $p$, (i) and (ii) are together equivalent to the condition
\begin{equation}\label{niceBLforTDI}
    \dim V \leq\frac{n+d}{d_h} \text{Avg}_j \dim B_{j, \theta} V, \forall \text{ subspace } V \subset \Bbb{R}^{d+n}.
\end{equation}

It has been proved that (\ref{niceBLforTDI}) is satisfied for all Parsell-Vinogradov manifolds, as long as one takes $M$ sufficiently large. As we have seen, this guarantees useful ball-inflation lemmas from the corresponding perturbed Brascamp-Lieb inequality. These lemmas then play key roles in proving sharp decoupling for these manifolds (Theorem 1.2 in \cite{guo2019integer}).

Finally we say a bit about the proof of (\ref{niceBLforTDI}) for these manifolds. When the dimension of the manifold $d=1$, this was done by Bourgain-Demeter-Guth \cite{bourgain2016proof} relatively easily via a Wronskian argument. Their method does not easily generalize to higher dimensions. The first interesting case there for $d=2, k=3$ is already quite different and was proved in the technical appendix of \cite{bourgain2017sharp} including a lot of case-by-case analysis.

In \cite{guo2019integer}, the full (\ref{niceBLforTDI}) for all $d$ and $k$ were verified in a systematic way. The statement boils down to lower bounding the rank of a kind of matrices whose entries are general $d$-variate polynomials. Here the matrices are viewed as ones over the field of rational functions in $d$ variables. The argument in \cite{guo2019integer} is involved, including algebraic considerations related to the method of degeneration and combinatorial tools related to the proof of the Schwartz-Zippel lemma. It was further extended in \cite{guo2020decoupling} to deal with more general TDI manifolds studied by Arkhipov-Chubarikov-Karatsuba and Parsell.

\subsection{A fruitful scale-dependent analysis} Although (\ref{niceBLforTDI}) holds for a great class of TDI manifolds, very unfortunately it also fails for many TDI manifolds. This may sound surprising at first sight: After all, we can allow as many caps as we want, and they only need to be in generic position! Nevertheless, one inevitably realizes the possible failure by considering the following example TDI manifold:
\begin{equation}\label{S1eq}
    S_1 = \{(x_1, x_2, x_3, x_1^2, x_2^2+x_1 x_3)\}.
\end{equation} Here $d=3$ and $n=2$. Take the order $h=1$ and $d_h = 3$.  A general tangent space is the span of three vectors 
\begin{equation}\label{vectorlist}
    (1, 0, 0, 2x_1, x_3), (0, 1, 0, 0, 2x_2) \text{ and }  (0, 0, 1, 0, x_1).
\end{equation}
Looking at for example the two dimensional $V = \{(*, 0, 0, *, 0)\}$, we see that 
$V$ is always orthogonal to the last two vectors in the list (\ref{vectorlist}). Hence $V$ has at most a one-dimensional projection to \emph{every} tangent space. From this we see (\ref{niceBLforTDI}) fails for $V$.

Now that (\ref{niceBLforTDI}) fails for many TDI manifolds like $S_1$, the ``classical'' Brascamp-Lieb inequality is not available to us anymore. But for all TDI manifolds including $S_1$, there are always sharp decoupling conjectures. What should we do to prove them? Fortunately, people recently found an alternative tool that is proved successful in a lot of situations where (\ref{niceBLforTDI}) fails. This tool is again a cousin of Brascamp-Lieb, known as regularized Brascamp-Lieb inequalities.

Recall from \S \ref{MKsec} and \S \ref{PBLsec} that in $\Bbb{R}^N$, for a Brascamp-Lieb datum $(B_1, \ldots, B_m, p_1, \ldots, p_m)$, one way to state the corresponding Brascamp-Lieb inequality is to set up families $\Bbb{T}_j$ containing $1$-neighborhoods of affine subspaces that are translations of $\ker B_j$ and quantify the overlapping of these families:

\begin{equation}\label{BLsetup}
    \int_{\Bbb{R}^N} \prod_{j=1}^m (\sum_{T \in \Bbb{T}_j} 1_T)^{p_j}  \leq BL(\textbf{B}, \textbf{p}) \prod_{j=1}^m |\Bbb{T}_j|^{p_j}.
\end{equation}

This has the disadvantage that one gets nothing at all if the Brascamp-Lieb constant $BL(\textbf{B}, \textbf{p})$ is infinity, which we see is already an issue when dealing with decoupling problems for manifolds such as $S_1$ in (\ref{S1eq}). On the other hand, in applications to Fourier restriction type problems such as decoupling, the ``tubes'' $T$ do not really have infinite lengths. Indeed, we usually care about the integral on the left hand side of (\ref{BLsetup}) inside some large ball $B_R$. This integral should be finite, regardless of whether $BL(\textbf{B}, \textbf{p}) < \infty$. It thus becomes meaningful to ask: If we replace the integral domain on the left hand side of (\ref{BLsetup}) by a ball $B_R$, can we come up with a similar right hand side with some constant depending on the radius $R$ of the ball even if $BL(\textbf{B}, \textbf{p})$ is infinite?

This problem is resolved by Maldague \cite{maldague2019regularized}, resulting in a theorem known as the \emph{regularized} Brascamp-Lieb inequality.

\begin{thm}[Regularized Brascamp-Lieb, \cite{maldague2019regularized}]\label{RBLthm}
In $\Bbb{R}^N$, for any Brascamp-Lieb datum $(B_1, \ldots, B_m, p_1, \ldots, p_m)$ and a given $R > 1$, the regularized Brascamp-Lieb constant is defined to be the smallest possible $C \geq 0$ such that
\begin{equation}\label{RBLineq}
    \int_{B_R^N} \prod_{j=1}^m f_j(B_j x)^{p_j}\leq C \prod_{j=1}^m (\int |f_j| )^{p_j}
\end{equation}
where each $f_j$ is constant on every lattice unit cube and a lattice unit cube in $\Bbb{R}^N$ is defined to be some $x+ [0, 1)^N$ where $x \in \Bbb{Z}^N$.

Then the regularized Brascamp-Lieb constant is $\sim_{(\textbf{B}, \textbf{p})} R^{\sup_{V \subset \Bbb{R}^n} \dim V -\sum_j p_j \dim (B_j V)}$.
\end{thm}

Theorem \ref{RBLthm} is sharp in terms of the power of $R$ one loses. Also from this theorem one can derive a version in terms of $1_T$ as in the setup above. We note that Maldague also proved the corresponding perturbed version with $R^\varepsilon$-loss in the same paper.

\begin{rem}
By rescaling, Theorem \ref{RBLthm} (or its perturbed counterpart) also implies versions where each $f_j$ is constant on lattice cubes of edge length $r<R$. When one takes $r\to 0$ and $R \to \infty$ this becomes Theorem \ref{BCCTfinite} concerning the original Brascamp-Lieb inequality. When $R \to \infty$ this is known as the discrete  Brascamp-Lieb inequality and when $r\to 0$ this  is known as the localized  Brascamp-Lieb inequality. Both latter cases were studied previously in \cite{bennett2005finite, bennett2008brascamp}.
\end{rem}

When Theorem \ref{RBLthm} first came out, it felt fundamental to people but did not see any immediate applications in Fourier analysis. Then in another few years, the authors of \cite{guo2021decoupling} realized that Theorem \ref{RBLthm} is exactly what is needed to compensate for the absence of ``classical'' Brascamp-Lieb in the proof of sharp decoupling for all quadratic TDI manifolds, including many ``degenerate'' ones like $S_1$ where (\ref{niceBLforTDI}) fails. Indeed, since one only needs ball-inflation lemmas to go from $\delta^{-1}$ balls to $\delta^{-2}$ balls, even if no useful Brascamp-Lieb is available, one can always look at the corresponding regularized Brascamp-Lieb inequalities and Theorem \ref{RBLthm} just gives the correct power of $\delta^{-1}$ one must lose. Would such an inequality be useful to produce sharp decoupling results like \cite{bourgain2015proof, bourgain2016proof}? The answer is affirmative: Guo, Oh, Zorin-Kranich and the author use Maldague's perturbed version of Theorem \ref{RBLthm} together with a carefully designed induction-on-scales machinery to prove the following very general theorem:

\begin{thm}[\cite{guo2021decoupling}]
For all $p, q \geq 2$, the sharp decoupling theorem (stated in details as Theorem 1.1 in \cite{guo2021decoupling}) holds for every TDI manifold of degree $2$.
\end{thm}

\cite{guo2021decoupling} settles the decoupling problem for all quadratic TDI manifolds (up to a possible $\delta^{-\varepsilon}$ loss), but their approach has not been generalized to higher degree cases. Looking forward, there are still lots of hope on the unexplored strength of regularized Brascamp-Lieb.

\subsection{Additional historical remarks} The general concept of TDI systems are introduced by Parsell-Prendiville-Wooley \cite{parsell2013near} where they explain the advantage of these systems and use the method of efficient congruencing to extend Vinogradov's result to higher dimensional  TDI systems. \cite{wooley2014translation} is a very good introduction to Vinogradov's work, related context of the circle method, TDI systems and the applications of efficient congruencing there. The ``intermediate dimension'' issue first showed up in the treatment of the cubic Parsell-Vinogradov system of Bourgain-Demeter-Guo in \cite{bourgain2017sharp}. For recent progress on decoupling for TDI manifolds, see \cite{bourgain2015proof, bourgain2016proof, bourgain2017decoupling, bourgain2016decouplings, bourgain2017decouplings, bourgain2017sharp, oh2018decouplings, guo2019integer, demeter2019sharp, guo2019decoupling, guo2020decoupling, guo2021decoupling}. A sharp decoupling inequality for degree $2$ TDI manifolds of dimension $3$ and codimension $2$ where  (\ref{niceBLforTDI}) fails like $S_1$ was obtained in \cite{guo2019decoupling} prior to the general result in \cite{guo2021decoupling}.

\section{The really sharp perturbed inequality}\label{reallysharpsec} In all discussions above, all perturbed versions of Brascamp-Lieb (such as Theorems \ref{MKthm} and \ref{PBLthm}) are stated with the left hand side on a ball of radius $R$ and they all have an $R^{\varepsilon}$-loss. This loss is not needed in the Brascamp-Lieb inequalities. Do we really need it in the perturbed version? One may suspect we do not since there seem no way to construct counterexamples. For instance, recall in examples in dimension two generalizing Figure \ref{MRpic}, one does not need to lose an $R^{\varepsilon}$ when estimating the overlap, which can be seen elementarily.

This is usually more of a mathematical curiosity to most people working in Fourier restriction theory, because as we have seen, very often one does not really mind these mild loss of $R^{\varepsilon}$ in applications. Nevertheless, people now know that we indeed do not need such a loss and we feature this story in this last section. The most influential part of it was the removal of the $R^{\varepsilon}$-loss in the multilinear Kakeya case by Guth \cite{guth2010endpoint}, where he took some of the most fundamental concepts in geometry and topology and applied them in analysis. This is an important step of the development of the \emph{polynomial method} in Fourier restriction. It is not possible to even give a sketch of that method here, so we will only say a bit about its earlier developments and how Guth introduced new ideas to prove the really sharp multilinear Kakeya (also known as endpoint multilinear Kakeya) inequality.

\subsection{The finite field Kakeya conjecture} The Kakeya conjecture asserts that as long as a compact set in $\Bbb{R}^n$ contains a unit line segment in every direction, that set must have full (Hausdorff) dimension $n$. If we look at the $\delta$-neighborhood of such a set, we see that this is a problem concerning how much translations of thin tubes in different directions can be compressed. Recall that we have discussed the importance of this kind of problems in Fourier restriction in the end of \S \ref{FRTprob}. More precisely, if we decompose $P^{n-1}$ into $R^{-\frac{1}{2}}$-caps, then due to the curvature of $P^{n-1}$, wave packets from different caps will be supported on tubes in different directions. Because of this connection, the Kakeya conjecture and a few related conjectures turn out to be the most important geometric measure theory questions cared by people working in Fourier restriction. See e.g. \cite{tao2001rotating} for an excellent introduction to Kakeya type problems and their connections to Fourier restriction.

The Kakeya conjecture is known in two dimensions \cite{davies1971some} but is widely open and notoriously difficult in dimensions $3$ and higher. Wolff \cite{wolff1999recent} brought up the analogue of the conjecture over finite fields. Initially that was again thought of as a difficult conjecture that may be well out of reach. But about 14 years ago Dvir surprised everyone by proving the conjecture over finite fields in a unique way.

\begin{thm}[Finite field Kakeya \cite{dvir2009size}]\label{FFF}
Fix the dimension $n$. For every finite field $\Bbb{F}$, if a subset of $\Bbb{F}^n$ contains a line in every direction, then it must have cardinality $\gtrsim_n |\Bbb{F}|^n$.
\end{thm}

Dvir's proof of Theorem \ref{FFF} is completely elementary, less than one page but very unusual to the people working on Fourier restriction and Kakeya. He constructs a nonzero ($n$-variate) polynomial that vanishes on the set. By counting degrees of freedom, if the set has small cardinality (denoted by $M$, much smaller than $|\Bbb{F}|^n$) then the polynomial can be chosen to have degree $O(M^{\frac{1}{n}})$, which is much lower than $|\Bbb{F}|$. Now for every direction the vanishing of the polynomial on a line in that direction implies the vanishing of the polynomial at the ``point at infinity''. But then such a polynomial must vanish at all points at infinity because of the Kakeya set property. This will violate the Schwartz-Zippel lemma and conclude the proof.

\begin{rem}
Construction of auxiliary polynomials had been useful in other problems in number theory and combinatorics. But this was the first time it was used on Kakeya or Fourier restriction type problems.
\end{rem}

People have then tried to adapt Dvir's proof to say things about the original Kakeya conjecture but it is still unclear as of today whether any new meaningful attack can be formed based on Dvir's proof alone. Nevertheless, his method did inspire a lot of new developments in the whole subject of Kakeya and Fourier restriction, starting with Guth's proof of the endpoint multilinear Kakeya conjecture.

\subsection{Endpoint multilinear Kakeya}\label{EMKsec} In \cite{guth2010endpoint}, Guth was able to remove the $R^{\varepsilon}$ loss in Theorem \ref{MKthm} and prove the following:

\begin{thm}[Endpoint multilinear Kakeya \cite{guth2010endpoint}]\label{EMKthm}
Fix the dimension $n$, then there exists a small $c(n)>0$ such that: Let $\Bbb{T}_j, 1\leq j \leq n,$ be a finite family of (infinite) unit tubes such that each $T \in \Bbb{T}_j$ has an angle $\leq c(n)$ against the $x_j$-direction. Then
\begin{equation}\label{EMKineq}
    \int_{\Bbb{R}^n} \prod_{j=1}^n (\sum_{T \in \Bbb{T}_j} 1_T)^{\frac{1}{n-1}} \lesssim \prod_{j=1}^n |\Bbb{T}_j|^{\frac{1}{n-1}}.
\end{equation}
\end{thm}

The proof of Guth is conceptually much related to Dvir's framework of the proof of finite field Kakeya, but it also relies on additional tools from geometry and algebraic topology. These tools themselves (e.g. cohomology of the projective space and the isoperimetric inequality) are known to be fundamental to people in geometry and topology but this was the first time they are used in Kakeya and Fourier restriction.

We sketch how to use Guth's idea to prove the following slightly easier conclusion than (\ref{EMKineq}): The size of the support of the function on the left hand side of (\ref{EMKineq}) is bounded by constant times the right hand side when each $|\Bbb{T}_j|$ are the same (denoted by $H$). In other words, we will explain how Guth's method can prove
\begin{equation}\label{toyEMKineq}
    \left|\bigcap_{j=1}^n \bigcup_{T \in \Bbb{T}_j} T\right| = O(H^{\frac{n}{n-1}})
\end{equation}
when each $|\Bbb{T}_j| = H$. It is good to compare this sketch with Dvir's proof.

Guth's approach to prove (\ref{toyEMKineq}) is as the following: Approximate the left hand side by a set of lattice unit cubes and we just need to bound the number of cubes (denoted by $M$) by the right hand side. By the \emph{polynomial ham-sandwich theorem} (that generalizes the usual ham-sandwich theorem and can be proved using cohomology of the real projective space with $\Bbb{Z}/2\mathbb{Z}$-coefficients), one can find a nonzero real polynomial $P$ of degree $O(M^{\frac{1}{n}})$ whose zero set \emph{bisects} each cube. By the \emph{isoperimetric inequality}, the zero set $Z(P)$ of the polynomial as a hypersurface must have volume $\gtrsim 1$ inside each cube. At each cube $Q$, we can assume $T_1, \ldots, T_n$ meet where $T_j \in \Bbb{T}_j$. Since $Z(P) \bigcap Q$ has volume $\gtrsim 1$ and the directions of $T_1, \ldots, T_n$ are quantitatively transverse, we see that $Z(P) \bigcap Q$ must have directed volume along the direction of some $T_j$ also $\gtrsim 1$. But $Z(P) \bigcap Q$ is contained in a thickening $\tilde{T}_j$ of $T_j$ by a constant times. Now choose a popular $j$ and we see when we sum over $T \subset \Bbb{T}_j$ of the directed volume of $Z(P) \bigcap \tilde{T}$ along the axis of  $\tilde{T}$, the result is $\gtrsim M$. But note that our polynomial $P$ can vanish at $\leq \deg P$ points on a generic line. By a simple integral geometric argument we see the contribution from each $T \in \Bbb{T}_j$ has to be $\leq \deg P = O(M^{\frac{1}{n}})$. Hence $M = O(M^{\frac{1}{n}}) H$, implying $M = O(H^{\frac{n}{n-1}})$.

To prove the full Theorem \ref{EMKthm}, Guth did extra technical work including introducing the concept of visibility, using the Lusternik-Schnirelmann vanishing lemma in algebraic topology, and applying ideas from geometry such as in Gromov's work \cite{gromov2003isoperimetry, gromov2010singularities}. Many main ideas in the full proof were already illustrated in the above toy model case though.

\cite{guth2010endpoint} has the significance of being the first example of using the zero set of real polynomials in Kakeya/Fourier restriction to obtain new and sharp results. Further developments of this idea have led to major breakthroughs such as the proof of the Erd\H{o}s distinct  distances conjecture \cite{guth2015erdHos} and the polynomial method in the Fourier restriction Conjecture (starting in \cite{guth2016restriction, guth2018restriction}).

We also remark that later Carbery-Valdimarsson gave a more elementary proof of Guth's endpoint multilinear Kakeya Theorem \ref{EMKthm} in \cite{carbery2013endpoint}. Their proof uses no algebraic topology beyond the Borsuk--Ulam theorem.

\subsection{Removing $R^{\varepsilon}$ in the full generality} For more general perturbed Brascamp-Lieb inequalities (Theorem \ref{PBLthm}), the $R^{\varepsilon}$ can be removed too. This was proved by the author \cite{zhang2017endpoint}:

\begin{thm}[Endpoint perturbed Brascamp-Lieb inequality \cite{zhang2017endpoint}]\label{EPBLthm}
Assumptions are the same as Theorem \ref{PBLthm}, on can replace $B_R$ by the whole $\Bbb{R}^n$, and remove the $R^{\varepsilon}$ factor in (\ref{PBLineq}).
\end{thm}

The proof idea of Theorem \ref{EPBLthm} was greatly inspired by Guth's proof of the endpoint multilinear Kakeya. We comment on one interesting new feature: As we see above, in Guth's proof one key component is to use the fact that generic lines can only intersect $Z(P)$ at $\leq \deg P$ points. The lines will be taken as translations of the core line in each tube $T$, and that the averaged number of intersection points with $Z(P)$ will then measure the directed volume of $Z(P) \bigcap T$. In the general Brascamp-Lieb setting, $T$ will be neighborhoods of higher dimensional subspaces (``core subspace'' of $T$, denoted by $c(T)$) and they  in general intersect $Z(P)$ at higher dimensional subvarieties and thus the Dvir-Guth framework needs modification. A natural thing one can do is to take more zero sets of polynomials $Z(P_1), \ldots, Z(P_m)$ where $m$ equals to the dimension of $c(T)$. Then generically the intersection of $Z(P_1), \ldots, Z(P_m)$ and a subspace parallel to $c(T)$ will be discrete sets of points that can then be counted, giving hope to extending Guth's approach. Finding $P_1, \ldots, P_m$ may seem a daunting task that could require new tools in algebraic geometry and topology. However, somewhat surprisingly, it turns out in \cite{zhang2017endpoint} that one can take $P_1, \ldots, P_m$ all to be generic perturbations of the same $P$ that was already constructed in the same way as in Guth's work \cite{guth2010endpoint}. By careful multilinear algebra considerations, one deduces Theorem \ref{EPBLthm} without the need of additional tools from geometry or topology.

Theorem \ref{EPBLthm} is stronger than Theorem \ref{PBLthm}, but we remark that its proof does rely on ingredients (the stability of the Brascamp-Lieb constant) that goes into the proof of Theorem \ref{PBLthm}.

Theorem \ref{EPBLthm} has also been generalized to the regularized setting (related to Theorem \ref{RBLthm}) by Zorin-Kranich \cite{zorin2019kakeya} via further developments of the ideas in \cite{guth2010endpoint, zhang2017endpoint} and the related \cite{carbery2013endpoint}. See \cite{gressman2022testing} for a related recent development.

\subsection{Additional historical remarks} The concepts of Kakeya sets and the Kakeya conjectures originated from \cite{kakeya1917some} more than a century ago. There is a related stronger Kakeya maximal function conjecture to the Kakeya set conjecture. It is known that the Fourier restriction Conjectures \ref{restrconj} and \ref{extconj} both imply the Kakeya conjecture. It was Bourgain \cite{bourga1991besicovitch} who first realized that progress on the Kakeya  maximal function conjecture can lead to progress on the Fourier restriction conjecture through the wave packet decomposition. These two problems are also closely connected to Bochner-Riesz and local smoothing. For example, the readers can use table 2 in \cite{tao1999bochner} to see known implications between them.

In Guth's paper \cite{guth2010endpoint}, as a warm-up he also first discussed how one can deduce the weaker conclusion (\ref{toyEMKineq}) and our presentation follows his. The polynomial ham-sandwich theorem was proved by Stone and Tukey in \cite{stone1942generalized}.

For the polynomial method, \cite{guth2016polynomial} is a good reference of it in combinatorics. Chapter 8 of \cite{demeter2020fourier} has a detailed account on how Guth applied this method to Fourier restriction in \cite{guth2016restriction}.

\section*{Acknowledgements} My research is supported by the NSF grant DMS-2207281. I would like to thank Rupert Frank and Shaoming Guo for comments and discussions. I am deeply indebted to the anonymous referee whose detailed suggestions and  historical remarks have greatly improved the accuracy and exposition of this article.

\newpage
\bibliographystyle{amsalpha}
\bibliography{BL}

\end{document}